\documentstyle[12pt]{article}

\newtheorem{defn0}{Definition}[section]
\newtheorem{prop0}[defn0]{Proposition}
\newtheorem{thm0}[defn0]{Theorem}
\newtheorem{lemma0}[defn0]{Lemma}
\newtheorem{corollary0}[defn0]{Corollary}
\newtheorem{example0}[defn0]{Example}
\newtheorem{conjecture0}[defn0]{Conjecture}

\newenvironment{defn}{\begin{defn0}}{\end{defn0}}
\newenvironment{prop}{\begin{prop0}}{\end{prop0}}
\newenvironment{thm}{\begin{thm0}}{\end{thm0}}
\newenvironment{lemma}{\begin{lemma0}}{\end{lemma0}}
\newenvironment{corollary}{\begin{corollary0}}{\end{corollary0}}
\newenvironment{example}{\begin{example0}\rm}{\end{example0}}

\newcommand{\defref}[1]{Definition~\ref{#1}}
\newcommand{\propref}[1]{Proposition~\ref{#1}}
\newcommand{\thmref}[1]{Theorem~\ref{#1}}
\newcommand{\lemref}[1]{Lemma~\ref{#1}}
\newcommand{\corref}[1]{Corollary~\ref{#1}}

\newcommand{\secref}[1]{Section~\ref{#1}}

\newcommand{\qed}{\mbox{~~~~\vrule height 1.2ex width .9ex depth .1ex}}

\newenvironment{proof}{\noindent {\bf Proof.}}{\qed\vskip 6pt}
\newenvironment{proofn}{\noindent {\bf Proof.}}{\vskip 6pt}

\newcommand{\mcly}{\mbox{\em Macaulay}}
\newcommand{\std}{Gr\"{o}bner}
\newcommand{\Std}{Gr\"{o}bner}

\newcommand{\mmbox}[1]{\mbox{${#1}$}}

\newcommand{\PP}{\mbox{\bf P}}

\newcommand{\ZZ}{\mbox{\bf Z}}

\newcommand{\proj}[1]{\mmbox{{\bf P}^{#1}}}
\newcommand{\affine}[1]{\mmbox{{\bf A}^{#1}}}

\newcommand{\sliver}{\hskip 0.01in}

\newcommand{\x}{\sliver{\bf x}\sliver}

\newcommand{\Oh}{{\cal O}}
\newcommand{\Ih}{{\cal I}}
\newcommand{\Fh}{{\cal F}}

\newcommand{\mult}{{\rm mult}}
\newcommand{\geomdeg}{\mbox{\rm geom-deg}}
\newcommand{\arithdeg}{\mbox{\rm arith-deg}}

\newcommand{\thuh}[1]{\raisebox{3pt}{\small {#1}}}

\newcommand{\p}{\sliver{\bf p}\sliver}
\newcommand{\q}{\sliver{\bf q}\sliver}

\newcommand{\ssubset}{\raisebox{-2.5pt}{$\,\stackrel\subset{\scriptstyle
\neq}\,$}}
\newcommand{\ssupset}{\raisebox{-2.5pt}{$\,\stackrel\supset{\scriptstyle
\neq}\,$}}

\newcommand{\quotes}[1]{\mbox{``#1''}}

\newcommand{\ini}{{\rm in}}
\newcommand{\sat}[1]{#1^{\rm sat}}
\newcommand{\setdef}[2]{\{\,#1\,|\,#2\,\}}

\newcommand{\Bigcap}[2]
{\raisebox{-5pt}{$\;\stackrel{\textstyle\stackrel{#2}\bigcap}{_{#1}}\;$}}

\newcommand{\Limits}[3]
{\raisebox{-5pt}{$\;\stackrel{\textstyle\stackrel{#3}{#1}}{_{#2}}\;$}}
\newcommand{\abs}[1]{\mmbox{\mid \! #1 \! \mid}}

\newcommand{\limt}{\lim_{t \rightarrow 0}}
\newcommand{\rz}{\!\!\mid_{t=0}}
\newcommand{\rx}{\!\!\mid_{t=1}}

\newcommand{\ord}{{\rm ord}}
\newcommand{\reg}{{\rm reg}}
\newcommand{\codim}{{\rm codim}}
\begin{document}

\title {What Can Be Computed \\
in Algebraic Geometry?}

\author {
Dave Bayer
\thanks {Partially supported by NSF grant DMS-90-06116.}
\and David Mumford
}

\date{May 7, 1992}

\maketitle
This paper evolved from a long series of discussions between the two
authors, going back to around 1980, on the problems of making
effective computations in algebraic geometry, and it took more definite
shape in a survey talk given by the second author at a conference on
Computer Algebra in 1984. The goal at that time was to bring together
the perspectives of theoretical computer scientists and of working
algebraic geometers, while laying out what we considered to be the
main computational problems and bounds on their complexity. Only
part of the talk was written down and since that time there has been a
good deal of progress. However, the material that was written up may
still serve as a useful introduction to some of the ideas and estimates
used in this field (at least the editors of this volume think so), even
though most of the results included here are either published
elsewhere, or exist as ``folk-theorems'' by now.

The article has four sections. The first two parts are concerned with
the theory of \std\ bases;  their construction provides the foundation
for most computations, and their complexity dominates the complexity
of most techniques in this area. The first part introduces \std\ bases
from a geometric point of view, relating them to a number of ideas
which we take up in more detail in subsequent sections. The second
part develops the theory of \std\ bases more carefully, from an
algebraic point of view. It could be read independently, and requires
less background. The third part is an investigation into bounds in
algebraic geometry of relevance to these computations. We focus on
the {\em regularity\,} of an algebraic variety (see \defref{def1}),
which, beyond its intrinsic interest to algebraic geometers, has
emerged as a measure of the complexity of computing \std\ bases (see
\cite{bs87a}, \cite{bs87b}, \cite{bs88}). A principal result in this part is
a bound on the regularity of any smooth variety by the second author:
\thmref{thm8}(b). This bound has stimulated subsequent work, and has now been
generalized by \cite{bel91}. Another result of interest is \propref{prop9},
which elucidates the scheme structure of the ideal membership problem. The
fourth part is a short discussion of work on algorithms for performing some
other key operations on varieties, some open problems about these operations
and some general ideas about what works and what doesn't, reflecting the
prejudices of the authors.

One of the difficulties in surveying this area of research is that
mathematicians from so many specialties have gotten involved, and
they tend both to publish in their own specialized journals and to have
specific agendas corresponding to their area. Thus one group of
researchers, the working algebraic geometers, are much more
interested in actually computing examples than in worst-case
complexity bounds. This group, including the first author, has put a
great deal of work into building a functioning system, \mcly, based on
\std\ bases, which has solved many problems and provided many
examples to the algebraic geometry community \cite{macaulay}.
Another group comes from theoretical computer science and is much
more interested in theoretical bounds than practical systems (cf. the
provocative comments in Lenstra's survey \cite{len92}). It seems to
us that more communication would be very helpful: On the one hand,
the working algebraic geometer knows lots of facts about varieties that
can be very relevant to finding fast algorithms. Conversely asymptotic
and/or worst-case performance bounds are sometimes,  at least,
important indicators of real-time performance. These  theoretical
bounds may also reveal important distinctions between classes of
procedures, and may pose new and deep problems in algebraic
geometry. Thus we will see in \secref{bounds} how regularity
estimates flesh out a picture explaining why \std\ basis computations
can have such explosive worst case behavior, yet be so useful for the
kinds of problems typically posed by mathematicians. Finally, to make
this article more useful in bridging this gap, we have tried to include a
substantial number of references in our discussions below.

\section{A Geometric Introduction}
\label{geometry}
Let $X$ be a subvariety or a subscheme of projective $n$-space $\proj
n$, over a field $k$. Let $\Fh$ be a vector bundle or a coherent sheaf
supported on $X$. We would like to be able to manipulate such objects
by computer. From algebra we get finite descriptions, amenable to
such manipulations: Let $S = k[x_0,\ldots, x_n]$ be the homogeneous
coordinate ring of $\proj n$. Then $X$ can be taken to be the
subscheme defined by a homogenous ideal $I \subset S$, and $\Fh$ can
be taken to be the sheaf associated to a finitely generated $S$-module
$M$. We can represent $I$ by a list of generators $(f_1, \ldots, f_r)$,
and $M$ by a presentation matrix $F$, where $$M_1
\stackrel{F}{\longrightarrow} M_0 \longrightarrow M \longrightarrow
0$$ presents $M$ as a quotient of finitely generated free $S$-modules
$M_0$, $M_1$. We concentrate on the case of an ideal $I$; by working
with the submodule $J = {\rm Im}(F) \subset M_0$, the module case
follows similarly.

The heart of most computations in this setting is a deformation of the
input data to simpler data, combinatorial in nature: We want to move
through a family of linear transformations of $\proj n$ so that in the
limit our objects are described by monomials. Via this family, we
hope to pull back as much information as possible to the original
objects of study.

Choose a one-parameter subgroup $\lambda(t) \subset GL(n+1)$ of the
diagonal form
$$\lambda(t) = \left[
\begin{array}{cccc}
t^{w_0} &  & &  \\
& t^{w_1} & &  \\
& & \ldots &  \\
& &  & t^{w_n}
\end{array}
\right],
$$ where $W = (w_0, \ldots, w_n)$ is a vector of integer weights. For
each $t \ne 0$, $\lambda(t)$ acts on $X$ via a linear change of
coordinates of $\proj n$, to yield the subscheme $X_t = \lambda(t) X
\cong X$. The limit $$X_0 = \limt X_t$$ is usually a simpler object,
preferable to $X$ for many computational purposes.

Even if we start out by restricting $X$ to be a subvariety rather than
a subscheme of $\proj n$, it does not suffice to take the limit $X_0$
set-theoretically; often all we will get pointwise in the limit is a
linear subspace $L \subset \proj n$, reflecting little besides the
dimension of the original variety $X$. By instead allowing this limit
to acquire embedded components and a nonreduced structure, we can
obtain an $X_0$ which reflects much more closely the character of $X$
itself.

We compute explicitly with the generators $f_1, \ldots, f_r$ of $I$:
Let $\lambda$ act on $S$ by mapping each $x_i$ to $t^{w_i} x_i$;
$\lambda$ maps each monomial $\x^A = x_0^{a_0} \cdots x_n^{a_n}$ to
$t^{W \cdot A} \x^A = t^{w_0 a_0 + \ldots + w_n a_n} x_0^{a_0} \cdots
x_n^{a_n}$. If $f = a \x^A + b \x^B + \ldots$, then $\lambda f = a
\,t^{W \cdot A} \x^A + b \,t^{W \cdot B} \x^B + \ldots$. We take the
projective limit $\ini(f) = \limt \lambda f$ by collecting the terms
of $\lambda f$ involving the least power of $t$; $\ini(f)$ is then the
sum of the terms $a \x^A$ of $f$ so $W \cdot A$ is minimal. For a
given $f$ and most choices of $\lambda$, $\ini(f)$ consists of a
single term.

The limit $X_0$ we want is defined with all its scheme structure by
the ideal $\ini(I) = \limt \lambda I$, generated by the set
$\setdef{\ini(f)}{f \in I}$. For a given $I$ and most choices of
$\lambda$, $\ini(I)$ is generated by monomials. Unfortunately, this
definition is computationally unworkable because $I$ is an infinite
set, and $\ini(I)$ need not equal $(\ini(f_1), \ldots, \ini(f_r))$ for
a given set of generators $f_1, \ldots, f_r$ of $I$. To understand how
to compute $\ini(I)$, we need to look more closely at the family of
schemes $X_t$ defined by $\lambda$.

Let $S[t]$ be the polynomial ring $k[x_0, \ldots, x_n, t]$; we view
$S[t]$ as the coordinate ring of a one-parameter family of projective
spaces ${\bf P}^n_t$ over the affine line with parameter $t$. For each
generator $f_j$ of $I$, rescale $\lambda f_j$ so the lowest power of
$t$ has exponent zero: Let $g_j = t^{-\ell} \lambda f_j$, where $\ell
= W \cdot A$ is the least exponent of $t$ in $\lambda f_j$. Then $f_j
= g_j\rx$ and $\ini(f_j) = g_j\rz$. Now, let $J \subset S[t]$ be the
ideal generated by $(g_1, \ldots, g_r)$; $J$ defines a family $Y$ over
$\affine 1$ whose central fiber is cut out by $(\ini(f_1), \ldots,
\ini(f_r))$.

What is wrong with the family $Y$? $Y$ can have extra components
over $t = 0$, which bear no relation to its limiting behavior as $t
\rightarrow 0$. Just as the set-theoretic limit $\limt X_t$ can be too
small (we need the nonreduced structure), this algebraically defined
limit can be too big; the natural limit lies somewhere in between.

The notion of a {\em flat\,} family captures exactly what we are
looking for here. For example, if $Y$ is flat, then there are no extra
components over $t = 0$. While the various technical definitions of
flatness can look daunting to the newcomer, intuitively flatness
captures exactly the idea that every fiber of a family is the natural
scheme-theoretic continuation of its neighboring fibers.

In our setting, all the $X_t$ are isomorphic for $t \ne 0$, so we only
need to consider flatness in a neighborhood of $t = 0$. Artin
\cite{art76} gives a criterion for flatness applicable here: The {\em
syzygies\,} of $g_1, \ldots, g_r$ are the relations $h_1 g_1 + \ldots
+ h_r g_r = 0$ for $h_1, \ldots, h_r \in S[t]$. Syzygies correspond to
elements $(h_1, \ldots, h_r)$ of the $S[t]$-module $S[t]^r$; the set
of all syzygies is a submodule of $S[t]^r$. $Y$ is a flat family at $t
= 0$ if and only if the restrictions $(h_1\rz, \ldots, h_r\rz)$ of
these syzygies to the central fiber generate the $S$-module of
syzygies of $g_1\rz, \ldots, g_r\rz$.

When $g_1\rz, \ldots, g_r\rz$ are single terms, their syzygies take on
a very simple form: The module of syzygies of two terms $a\x^A$,
$b\x^B$ is generated by the syzygy $b\x^C (a\x^A) - a\x^D (b\x^B) = 0$,
where $\x^E = \x^C \x^A = \x^D \x^B$ is the least common multiple of
$\x^A$ and $\x^B$. The module of syzygies of $r$ such terms is
generated (usually not minimally) by the syzygies on all such pairs.

We want to lift these syzygies to syzygies of $g_1, \ldots, g_r$,
working modulo increasing powers of $t$ until each syzygy lifts
completely.  Whenever we get stuck, we will find ourselves staring at
a new polynomial $g_{r+1}$ so $t^\ell g_{r+1} \in J$ for some $\ell >
0$. Including $g_{r+1}$ in the definition of a new $J^\prime \supset
J$ has no effect on the family defined away from $t = 0$, but will cut
away unwanted portions of the central fiber; what we are doing is
removing $t$-torsion.  By iterating this process until every syzygy
lifts, we obtain explicit generators $g_1, \ldots, g_r, g_{r+1},
\ldots, g_s$ for a flat family describing the degeneration of $X =
X_1$ to a good central fiber $X_0$.  The corresponding generators
$g_1\rx, \ldots, g_s\rx$ of $I$ are known as a {\em \std\ basis} for
$I$.

This process is best illustrated by an example. Let $S = k[w,x,y,z]$
be the coordinate ring of $\proj 3$, and let $I = (f_1, f_2, f_3)
\subset S$ for $$f_1= w^2-xy,\; f_2 = wy-xz,\; f_3 = wz-y^2.$$ $I$
defines a twisted cubic curve $X \subset \proj 3$; $X$ is the image of
the map $(r,s) \mapsto (r^2s,r^3,rs^2,s^3)$. Let
$$\lambda(t) = \left[
\begin{array}{cccc}
t^{-16} &  & &  \\
& t^{-4} & &  \\
&  & t^{-1} &  \\
& &  & t^0
\end{array}
\right].
$$
If $w^a x^b y^c z^d$ is a monomial of degree $< 4$, then $\lambda
\cdot w^a x^b y^c z^d = t^{-\ell} w^a x^b y^c z^d$ where $\ell =
16a + 4b + c$. Thus, sorting the monomials of $S$ of each degree $<
4$ by increasing powers of $t$ with respect to the action of $\lambda$
is equivalent to sorting the monomials of each degree in lexicographic
order.

We have
\begin{eqnarray*}
g_1 & =\; t^{32}\lambda f_1 & =\; w^2-t^{27}xy, \\
g_2 & =\; t^{17}\lambda f_2 & =\; wy-t^{13}xz, \\
g_3 & =\; t^{16}\lambda f_3 & =\; wz-t^{14}y^2.
\end{eqnarray*}
The module of syzygies on $w^2$, $wy$, $wz$ is generated by the three
possible pairwise syzygies; we start with the syzygy $y(w^2) - w(wy) =
0$. Substituting $g_1$, $g_2$ for the lead terms $w^2$, $wy$ we get
$$y(w^2-t^{27}xy) - w(wy-t^{13}xz) = t^{13}wxz - t^{27}xy^2$$ which is
a multiple $t^{13}x$ of $g_3$. Thus, the syzygy $$y g_1 - w g_2 -
t^{13} x g_3 = 0$$ of $g_1$, $g_2$, $g_3$ restricts to the monomial
syzygy $y(w^2) -w(wy) = 0$ when we substitute $t = 0$, as desired.

Similarly, the syzygy $$z g_1 - t^{14}y g_2 - w g_3 = 0$$ restricts to
the monomial syzygy $z(w^2) -w(wz) = 0$. When we attempt to lift
$z(wy) - y(wz) = 0$, however, we find that $$z(wy-t^{13}xz) -
y(wz-t^{14}y^2) = - t^{13}xz^2 + t^{14}y^3.$$ $xz^2$ is not a multiple
of $w^2$, $wy$, or $wz$, so we cannot continue; $J = (g_1, g_2, g_3)$
does not define a flat family. Setting $t = 1$, the troublesome
remainder is $- xz^2 + y^3$. Making this monic, let $f_4 = xz^2 -
y^3$; $f_4 \in I$ and $$g_4 =\; t^{4}\lambda f_4 =\; xz^2 - ty^3.$$
Adjoin $g_4$ to the ideal $J$, redefining the family $Y$. Now, $$z g_2
- y g_3 + t^{13} g_4 = 0$$ restricts to $z(wy) - y(wz) = 0$ as
desired.

The module of syzygies of $w^2$, $wy$, $wz$, and $xz^2$ is generated
by the pairwise syzygies we have already considered, and by the syzygy
$xz(wz) - w(xz^2) = 0$, which is the restriction of $$-ty^2 g_2 +xz
g_3 - w g_4 = 0.$$ Thus, $J = (g_1, g_2, g_3, g_4)$ defines a flat
family $Y$, and $$w^2-xy, \;wy-xz, \;wz-y^2, \;xz^2-y^3$$ is a \std\
basis for $I$. The limit $X_0$ is cut out by the monomial ideal
$\ini(I) = (w^2, wy, wz, xz^2)$, which we shall see shares many
properties with the original ideal $I$. Note that $xz^2 - y^3 = 0$
defines the projection of $X$ to the plane $\proj 2$ in $x$, $y$, and
$z$.

The scheme structure of $X_0$ is closely related to the combinatorial
structure of the monomial $k$-basis for $S/\ini(I)$: For each degree
$d$ in our example, the monomials not belonging to $\ini(I)$ consist
of three sets $\{x^d, x^{d-1}y,\ldots,y^d\}$, $\{x^{d-1}z,
x^{d-2}yz,\ldots,y^{d-1}z\}$, $\{y^d, y^{d-1}z,\ldots,z^d\}$, and a
lone extra monomial $x^{d-1}w$. The first two sets correspond to a
double line supported on $w = z = 0$, the third set to the line $w = x
= 0$, and the extra monomial to an embedded point supported at $w = y
= z = 0$.  Together, this describes the scheme structure of $X_0$. The
first two sets consist of $d+1$ and $d$ monomials, respectively; the
third set adds $d-1$ new monomials, and overlaps two monomials we have
already seen. With the extra monomial, we count $3d+1$ monomials in
each degree, which agrees with the dimensions of the graded pieces of
$S/I$.  The embedded point is crucial; it makes this count come out
right, and it alone keeps $X_0$ nonplanar like $X$.

The new monomial generator $xz^2$ of $\ini(I)$ excludes the line $w =
y = 0$ from $X_0$; combinatorially, it excludes all but three
monomials of the set $\{x^d, x^{d-1}z,\ldots,z^d\}$ from the
monomial $k$-basis for each degree of the quotient $S/\ini(I)$. We can
see that this line is unwanted as follows: Away from $t = 0$, $Y$ is
parametrized by $(r,s,t) \mapsto (t^{16}r^2s, t^4r^3,trs^2,s^3,t)$.
Thus, fixing $r$ and $s$, the curve $(r,ts,t) \mapsto (t^{17}r^2s,
t^4r^3,t^3rs^2,t^3s^3,t)$, with projective limit $(0,0,r,s,0)$ as $t
\rightarrow 0$. Similarly, the curve $(r,t^3s,t^2)$ has as its limit
$(0,r^2,s^2,0,0)$. These calculations show that the lines $w = z = 0$ and
$w = x = 0$ indeed belong set-theoretically to the limit $X_0$. We can
find no such curve whose limit is a general point on the line $w = y =
0$, for $(r, t^4s, t^3)$ doesn't work. Thus, the line $w = y = 0$ sticks out
of the good total space $Y$.

One usually computes \std\ bases by working directly in the ring $S$,
dispensing with the parameter $t$. The one-parameter subgroup
$\lambda$ is replaced by a total order on the monomials of each
degree, satisfying the {\em multiplicative\,} property $\x^A > \x^B
\Rightarrow \x^C \x^A > \x^C \x^B$ for all $\x^C$.  In fact, for our
purposes these are equivalent concepts: The weight vector $W$
associated with $\lambda$ induces the order $\x^A > \x^B
\Longleftrightarrow W\cdot A < W\cdot B$, which is a total
 multiplicative order in low degrees as long as no two monomials have
the same weight.  Conversely, given any multiplicative order and a
degree bound $d$, one can find many $\lambda$ which induce this order
on all monomials of degree $< d$. See \cite{bay82}, \cite{rob85} for
characterizations of such orders.

We shall be particularly interested in two multiplicative orders, the
{\em lexicographic\,} order used in our example, and the {\em
reverse lexicographic\,} order. The lexicographic order simply expands
out the monomials of each degree into words, and sorts them
alphabetically, i.e. $\x^A > \x^B$ iff the first nonzero entry in $A-B$ is
positive. The reverse lexicographic order pushes highest
powers of $x_n$ in any expression back to the end, then within these
groups pushes highest powers of $x_{n-1}$ to the end, etc., i.e.
$\x^A > \x^B$ iff the last nonzero entry of $A-B$ is negative.

What do these orders mean geometrically? The dominant effect of the
lexicographic order is a projection from $\proj n$ to $\proj{n-1}$,
eliminating $x_0$.  A second order effect is a projection to
$\proj{n-2}$, and so forth. We could compute the deformation from $X$
to $X_0$ with respect to the lexicographic order in stages carrying
out these projections, first applying a $\lambda$ with $W =
(-1,0,\ldots,0)$, then with $W = (-1,-1,0,\ldots,0)$, etc.
Alternatively, for monomials of each degree $< d$, we can apply the
single $\lambda$ with $W = (-d^{n-1},\ldots,-d,-1,0)$, generalizing
the $\lambda$ used in our example. Use of the lexicographic order
tends to muck up the family $Y$ more than necessary in most
applications, because projections tend to complicate varieties.

For the reverse lexicographic order, the dominant effect is a
projection of $\proj n$ down to the last coordinate point
$(0,\ldots,0,1)$. As a secondary effect, this order projects down to
the last coordinate line, and so forth. In other words, this order
first tries to make $X$ into a cone over the last coordinate point,
and only then tries to squash the result down to or cone it over the
last coordinate line, etc.  For monomials of each degree $< d$, this
can be realized by applying $\lambda$ with $W =
(0,1,d,\ldots,d^{n-1})$. Like such cones, the reverse lexicographic
order enjoys special properties with respect to taking linear sections
of $X$ or $X_0$ by intersection with the spaces defined by the last
variable(s) (see \cite{bs87a}). The preferred status of the reverse
lexicographic order can be attributed to this relationship, because
generic linear sections do not complicate varieties.

For example, if we take $X$ to be three general points in $\proj 2$, then using
the lexicographic order $X_0$ becomes a triple point on a line, because the
first order effect is the projection of the three points to a line, and the
second order limiting process keeps the points within this line. By contrast,
if we use the reverse lexicographic order then $X_0$ becomes the complete first
order neighborhood of a point (a point doubled in all directions). This is
because the first order limiting process brings the three points together from
distinct directions, tracing out a cone over the three points. The first order
neighborhood of the vertex in this cone has multiplicity 3, and is the same as
the complete first order neighborhood in the plane of this vertex.

For those familiar with the theory of valuations in birational geometry
\cite[Vol. II, Ch. VI]{zs76}, the lexicographic and reverse lexicographic
orders have simple interpretations. Recall that if $X$ is a variety of
dimension $n$, and $$F: X = Z_0 \supset Z_1 \supset Z_2 \supset \ldots
\supset Z_n$$ is a flag of subvarieties, $\codim_X(Z_i) = i$, with
$Z_i$ smooth at the generic point of $Z_{i+1}$, then we can define a
rank $n$ valuation $v_F$ on $X$ as follows: For each $i = 1,\ldots, n-1$,
fix $f_i$ to be a function on $Z_{i-1}$ with a $1\thuh{st}$ order zero on
$Z_i$.  Then for any function $f$, we can define $e_1 = \ord_{Z_1}(f)$,
$e_2 =
\ord_{Z_2}((f/f_1^{e_1})\!\mid_{Z_1})$, etc.,
and $v_F(f) = (e_1, \ldots, e_n) \in \ZZ^n$, where the value group
$\ZZ^n$ is ordered lexicographically. The arbitrarily chosen $f_i$ are
not needed to compare two functions $f$, $g$: We have $v_F(f) \succ
v_F(g)$ if and only if $\ord_{Z_1}(f/g) > 0$, or if this order is zero and
$\ord_{Z_2}((f/g)\!\mid_{Z_1}) > 0$, and so forth. Such a valuation also
defines an order on each graded piece $S_d$ of the homogeneous
coordinate ring: take any $f_0 \in S_d$ and say $f > g$ if and only if
$v_F(f/f_0) \succ v_F(g/f_0)$. More generally, one may take the $Z_i$
to be subvarieties of a variety $X^\prime$  dominating $X$ and pull
back functions to $X^\prime$ before computing $v_F$.

The lexicographic order on monomials of each degree of $\proj n$ is
now induced by the flag $$\proj n \supset V(x_0) \supset V(x_0,x_1)
\supset \ldots \supset V(x_0,\ldots,x_{n-1}).$$ For example, the first
step in the comparison defining $v_F(\x^A/f_0) \succ v_F(\x^B/f_0)$
has the effect of asking if $a_0 - b_0 > 0$.

The reverse lexicographic order is induced by a flag on a blowup $X$
of $\proj n$: First blow up $V(x_0,\ldots,x_{n-1})$ and let $E_1$ be the
exceptional divisor. Next blow up the proper transform
of $V(x_0,\ldots,x_{n-2})$, and let $E_2$ be this exceptional
divisor. Iterating, we can define a flag
$$X \supset E_1 \supset E_1 \cap E_2 \supset \ldots \supset E_1 \cap
\ldots \cap E_n$$
which induces the reverse lexicographic order on monomials in each
degree.  For example, looking at the affine piece of the first blow up
obtained by substituting $x_0 = x_0^\prime x_{n-1}, \;\ldots, \;x_{n-2}
= x_{n-2}^\prime x_{n-1}$, the power of $x_{n-1}$ in the transform of
$\x^A$ is $a_0 + \ldots + a_{n-1}$, which is the order of vanishing of
this monomial on $E_1$. Thus, the first step in the comparison defining
$v_F(\x^A/f_0) \succ v_F(\x^B/f_0)$ has the effect of asking if $a_0 +
\ldots + a_{n-1} - b_0 - \ldots - b_{n-1} > 0$, which is what we want.

Taking into account the equivalence between multiplicative orders and
one-parameter subgroups, the process we have described in $S[t]$ is
exactly the usual algorithm for computing \std\ bases. It is
computationally advantageous to set $t = 1$ and dismiss our extra
structure as unnecessary scaffolding, but it is conceptually
advantageous to treat our viewpoint as what is ``really'' going on;
many techniques of algebraic geometry become applicable to the family
$Y$, and assist in analyzing the complexity of \std\ bases.
Moreover, this picture may help guide improvements to
the basic algorithm. For example, for very large problems, it could be
computationally more efficient  to degenerate to $X_0$
in several stages; this has not been tried in practice.

The coarsest measure of the complexity of a \std\ basis
is its maximum degree, which is the highest degree of a generator of
the ideal $\ini(I)$ defining $X_0$. This quantity is bounded by the
better-behaved {\em regularity\,} of $\ini(I)$: The regularity of an
ideal $I$ is the maximum over all $i$ of the degree minus $i$ of any
minimal $i$\thuh{th} syzygy of $I$, treating generators as $0$\thuh{th}
syzygies. When $I$ is the largest (the {\em saturated}) ideal defining
a scheme $X$, we call this the regularity of $X$. We take up
regularity in detail in \secref{bounds}; here it suffices to know that
regularity is {\em upper semi-continuous\,} on flat families, i.e. the
regularity can only stay the same or go up at special fibers.

Let $\reg(I)$ denote the regularity of $I$, and $\reg_0(I)$ denote the
highest degree of a generator of $I$.  In our case, $t = 0$ is the
only special fiber, and the above says that $$\reg_0(I) \;\le\;
\reg(I) \;\le\; \reg(\ini(I)) \;\ge\; \reg_0(\ini(I)),$$ where
$\reg_0(I)$ can be immediately determined from the input data, and
$\reg_0(\ini(I))$ is the degree-complexity of the \std\ basis
computation.  In practice, each of these inequalities are often
strict.

However when $k$ is infinite, then for any set of coordinates for
$\proj n$ chosen from a dense open set $U \subset GL(n+1)$ of
possibilities, Galligo (\cite{gal74}; see also \cite{bs87b}) has shown
that the limiting ideal $\ini(I)$ takes on a very special form:
$\ini(I)$ is invariant under the action of the Borel subgroup of upper
triangular matrices in $GL(n+1)$.  This imposes strong geometric
conditions on $X_0$. In particular, the associated primes of $\ini(I)$
are also Borel-fixed, so they are all of the form $(x_0, \ldots, x_i)$
for various $i$. This means that the components of $X_0$ are supported
on members of a flag.

In characteristic zero, it is shown in \cite{bs87a} that the
regularity of a Borel-fixed ideal is exactly the maximum of the
degrees of its generators, or in our notation, that $\reg(\ini(I)) =
\reg_0(\ini(I))$ when $\ini(I)$ is Borel-fixed. Thus, for generic
coordinates in characteristic zero, the degree-complexity of computing
\std\ bases breaks down into two effects: the gap $\reg_0(I) \le
\reg(I)$ between the input degrees and the regularity of $X$, and the
gap $\reg(I) \le \reg(\ini(I))$ allowed by upper-semicontinuity.

A combination of theoretical results, hunches and experience guides the
practitioner in assessing the first gap; what about the second? Does the
regularity have to jump at all? One can easily find examples of ideals and
total orders exhibiting such a jump, but in \cite{bs87a}, it is shown that for
the reverse lexicographic order, in generic coordinates and any characteristic,
there is no jump: $\reg(I) = \reg(\ini(I))$, so in characteristic zero we have
$$\reg_0(\ini(I)) = \reg(I).$$ In this sense, this order is an optimal choice:
{\em For the reverse lexicographic order, the degree-complexity of a \std\
basis computation is exactly the regularity of the input data.} This agrees
with experience; computations made on the same inputs using the lexicographic
order can climb to much higher degrees than the reverse
lexicographic order, in practice.

For many applications, one is free to choose any order, but some
problems restrict us to using orders satisfying combinatorial
properties which the reverse lexicographic order fails to satisfy. An
example, developed further in \secref{division}, is that of
eliminating variables, or equivalently, of computing projections. To
compute the intersection of $I$ with a subring $R =
k[x_i,\ldots,x_n]$, it is necessary to use an order which in each
degree sorts all monomials not in $R$ ahead of any monomial in $R$.
The lexicographic order is an example of such an order, for each $i$
simultaneously. This strength comes at a cost; we are paying in
regularity gaps for properties we may not need in a particular
problem. An optimal order if you need one specific projection (in the
same sense as above) is constructed by sorting monomials by total
degree in the variables to be eliminated, and then breaking ties using
the reverse lexicographic order. See \cite{bs87b} for this result, and
a generalization to the problem of optimally refining any nonstrict
order.

Using this elimination order, one finds that the inherent
degree-complexity of a computation is given not by the regularity of
$X$ itself, but rather by the regularity of the {\em flat
projection\,} $X^\prime$ of $X$, which is the central fiber of a flat
family which animates the desired projection of $X$ as $t \rightarrow
0$. The jump in regularity between $X$ and $X^\prime$ is unavoidable;
by choosing an optimal order, we avoid the penalty of a further jump
in regularity between $X^\prime$ and $X_0$.

The regularity of algebraic varieties or schemes $X$ is far from being
well understood, but there is considerable interest in its study; this
computational interpretation of regularity as the inherent
degree-complexity of an ideal is but one more log on the fire.

{}From a theoretical computer science perspective, the full complexity of
computing \std\ bases is determined not merely by the highest degree
$\reg_0(I)$ in the basis, but by the total number of arithmetic operations in
the field $k$ required to compute this basis. This has not been analyzed in
general, but for $0$-dimensional ideals $I$, Lakshman and Lazard (\cite{lak91},
\cite{ll91}) have shown that the complexity of computing reduced \std\ bases is
bounded by a polynomial in $d^n$, where $d$ is the maximum degree of the
generators, and $n$ is the number of variables.


\section{\Std\ Bases}
\label{division}
Let $S = k[x_0, \ldots, x_n]$ be a graded polynomial ring over the
field $k$, and let $I \subset S$ be a homogeneous ideal.

Let $S_d$ denote the finite vector space of all homogeneous, degree d
polynomials in $S$, so $S = S_0 \oplus S_1 \oplus \ldots \oplus S_d
\oplus \ldots $. Writing $I$ in the same manner as $I = I_0 \oplus I_1 \oplus
\ldots \oplus I_d \oplus \ldots $, we have $I_d \subset S_d$ for each $d$.
Recall that the Hilbert function of $I$ is defined to be the function
$p(d) = \dim(I_d)$, for $d \ge 0$.

A total order $>$ on the monomials of $S$ is said to be {\it
multiplicative\,} if whenever $\x^A > \x^B$ for two monomials $\x^A$,
$\x^B$, then $\x^C \x^A > \x^C \x^B$ for all monomials $\x^C$. This
condition insures that if the terms of a polynomial are in order with
respect to $>$, then they remain in order after multiplication by a
monomial.

\begin{defn}
\label{id1}
Let $>$ be a multiplicative order. For a homogeneous polynomial $f =
c_1 \x^{A_1} + \ldots + c_m \x^{A_m}$ with $\x^{A_1} > \ldots >
\x^{A_m}$, define the initial term $\ini(f)$ to be the lead (that is, the
largest) term $c_1 \x^{A_1}$ of $f$.  For a homogeneous ideal $I
\subset S$, define the initial ideal $\ini(I)$ to be the monomial ideal
generated by the lead terms of all elements of $I$.
\end{defn}

Note that the definitions of $\ini(f)$ and $\ini(I)$ depend on the
choice of multiplicative order $>$. See \cite{bm88} and \cite{mr88}
for characterizations of the finite set of $\ini(I)$ realized as the
order $>$ varies.

Fix a multiplicative order $>$ on $S$.

\begin{prop}[Macaulay]
\label{id2}
I and \ini(I) have the same Hilbert function.
\end{prop}

\begin{proof}
(\cite{mac27}) The lead terms of $I_d$ span $\ini(I)_d$, because every
monomial $\x^A \in \ini(I)$ is itself the lead term $\ini(f)$ of some
polynomial $f \in I$: Since $\x^A = \x^C \x^B$ for some $\x^B =
\ini(g)$ with $g \in I$, we have $\x^A = \ini(f)$ for $f = \x^C g$.

Choose a $k$-basis $B_d \subset I_d$ with distinct lead terms, and let
$\ini(B_d)$ be the set of lead terms of $B_d$; $\ini(B_d)$ has
cardinality $p(d) = \dim(I_d)$. Since any element of $I_d$ is a linear
combination of elements of $B_d$, any lead term of $I_d$ is a scalar
multiple of an element of $\ini(B_d)$. Thus, $\ini(B_d)$ is a basis
for $\ini(I)_d$, so $p(d) = \dim(\ini(I)_d).$
\end{proof}

One can compute the Hilbert function of $I$ by finding $\ini(I)$ and
applying this result; see \cite{mm83}, \cite{bcr91}, and \cite{bs92}.

\begin{corollary}
\label{id3}
The monomials of S which don't belong to $\ini(I)$ form a $k$-basis
for $S /I$.
\end{corollary}

\begin{proof}
These monomials are linearly independent in $S / I$, because any
linear relation among them is a polynomial belonging to $I$, and all
such polynomials have lead terms belonging to $\ini(I)$. These
monomials can be seen to span $S / I$ by a dimension count, applying
\propref{id2}.
\end{proof}

Two examples of multiplicative orders are the lexicographic order and
the reverse lexicographic order. $\x^A > \x^B$ in the lexicographic
order if the first nonzero coordinate of $A-B$ is positive. For
example, if $S = k[w, x, y, z]$, then $w > x > y > z$ in $S_1$, and
$$w^2 > wx > wy > wz > x^2 > xy > xz > y^2 > yz > z^2$$ in $S^2$.

$\x^A > \x^B$ in the reverse lexicographic order if the last nonzero
coordinate of $A-B$ is negative. For example, if $S = k[w, x, y, z]$,
then $w > x > y > z$ in $S_1$, and $$w^2 > wx > x^2 > wy > xy > y^2 >
wz > xz > yz > z^2$$ in $S^2$. These two orders agree on $S_1$, but
differ on the monomials of $S$ of degree $> 1$ when $n \ge 2$.

The lexicographic order has the property that for each subring $k[x_i,
\ldots, x_n] \subset S$ and each polynomial $f \in S$, $f \in k[x_i,
\ldots, x_n]$ if and only if $\ini(f) \in k[x_i, \ldots, x_n]$. The
reverse lexicographic order has the property that for each $f \in
k[x_0, \ldots, x_i]$, $x_i$ divides $f$ if and only if $x_i$ divides
$\ini(f)$.

One can anticipate the applications of these properties by considering
a $k$-basis $B_d \subset I_d$ with distinct lead terms, as in the
proof of \propref{id2}. With respect to the lexicographic order, $B_d
\cap k[x_i, \ldots, x_n]$ is then a $k$-basis for $I_d \cap k[x_i,
\ldots, x_n]$ for each $i$. With respect to the reverse lexicographic
order, $B_d \cap (x_n)$ is then a $k$-basis for $I_d \cap (x_n)$.
Thus, these orders enable us to find polynomials in an ideal which do
not involve certain variables, or which are divisible by a certain
variable. For a given degree $d$, one could construct such a basis
$B_d$ by applying Gaussian elimination to an arbitrary $k$-basis for
$I_d$. However, this cannot be done for all $d$ at once; such a
computation would be infinite. We will finesse this difficulty by
instead constructing a finite set of elements of I whose monomial
multiples yield polynomials in $I$ with every possible lead term.

Such sets can be described as follows:

\begin{defn}
\label{id4}
A list $F = [f_1, \ldots, f_r] \subset I$ is a (minimal) \std\ basis
for $I$ if $\ini(f_1), \ldots, \ini(f_r)$ (minimally) generate
$\ini(I)$.
\end{defn}

$\ini(I)$ is finitely generated because $S$ is Noetherian, so \std\
bases exist for any ideal I.

The order of the elements of $F$ is immaterial to this definition, so
$F$ can be thought of as a set. We are using list notation for $F$
because we are going to consider algorithms for which the order of the
elements is significant. For convenience, we shall extend the notation
of set intersections and containments to the lists $F$.

A minimal set of generators for an ideal $I$ need not form a \std\
basis for I. For example, if $S = k[x, y]$ and $I = (x^2 + y^2, xy)$,
then with respect to the lexicographic order, $\ini(x^2 + y^2) = x^2$
and $\ini(xy) = xy$. Yet $y(x^2 + y^2) - x(xy) = y^3 \in I$, so $y^3
\in \ini(I)$. Thus, any \std\ basis for $I$ must include $y^3$; it can
be shown that $\ini(I) = (x^2, xy, y^3)$ and $[x^2 + y^2, xy, y^3]$ is
a \std\ basis for $I$.

On the other hand,

\begin{lemma}
\label{id5}
If $F = [f_1, \ldots, f_r]$ is a \std\ basis for $I$, then $f_1,
\ldots, f_r$ generate $I$.
\end{lemma}

\begin{proof}
For each degree $d$, we can construct a $k$-basis $B_d \subset I_d$
with distinct lead terms, whose elements are monomial multiples of
$f_1, \ldots, f_r$: For each $\x^A \in \ini(I)_d$, $\x^A$ is a scalar
multiple of $\x^C \ini(f_i)$ for some $\x^C$ and some $i$; include
$\x^C f_i$ in the set $B_d$. Thus, the monomial multiples of $f_1,
\ldots, f_r$ span $I$.
\end{proof}

\begin{prop}[Spear, Trinks]
\label{id6}
Let $R \subset S$ be the subring $R = k[x_i, \ldots, x_n]$. If $F =
[f_1, \ldots, f_r ]$ is a \std\ basis for the ideal $I$ with respect
to the lexicographic order, then $F \cap R$ is a \std\ basis for the
ideal $I \cap R$. In particular, $F \cap R$ generates $I \cap R$.
\end{prop}

\begin{proof}
(\cite{spe77}, \cite{zac78}, \cite{tri78}) Let $f \in I \cap R$;
$\ini(f)$ is a multiple of $\ini(f_i)$ for some $i$. Since $\ini(f)
\in R$, $\ini(f_i) \in R$, so $f_i \in R$. Thus, $F \cap R$ is a \std\
basis for $I \cap R$. By \lemref{id5}, $F \cap R$ generates $I \cap
R$.
\end{proof}

\propref{id6} has the following geometric application: If $I$ defines
the subscheme $X \subset \PP^n$, then $I \cap k[x_i, \ldots, x_n]$
defines the projection of $X$ to $\PP^{n-i} = {\rm Proj}(k[x_i,
\ldots, x_n])$.

Recall that the saturation $\sat I$ of $I$ is defined to be the
largest ideal defining the same subscheme $X \subset \PP^n$ as $I$.
$\sat I$ can be obtained by taking an irredundant primary
decomposition for $I$, and removing the primary ideal whose associated
prime is the irrelevant ideal $(x_0, \ldots, x_n)$. $I$ is saturated
if $I = \sat I$.

If the ideal $I$ is saturated, and defines a finite set of points $X
\subset \PP^n$, then $I \cap k[x_{n-1}, x_n]$ is a principal ideal
$(f)$, where $\{f=0\}$ is the image of the projection of $X$ to $\PP^1
= {\rm Proj}(k[x_{n-1}, x_n])$. Given a linear factor of $f$ of the
form $(b x_{n-1} - a x_n)$, we can make the substitution $x_{n-1} = a
z$, $x_n = b z$ for a new variable $z$, to obtain from $I$ an ideal $J
\subset k[x_0, \ldots, x_{n-2}, z]$ defining a finite set of points in
$\PP^{n-1}$. For each point $(c_0, \ldots, c_{n-2}, d)$ in the zero
locus of $J$, $(c_0, \ldots, c_{n-2}, ad, bd)$ is a point in the zero
locus of $I$.

If $X \subset \PP^{n-1}$ is of dimension $1$ or greater, then in
general $I \cap k[x_{n-1}, x_n] = (0)$, because a generic projection
of $X$ to $\PP^1$ is surjective. In this case, an arbitrary
substitution $x_{n-1} = a z$, $x_n = b z$ can be made, and the process
of projecting to $\PP^1$ iterated. Thus, the lexicographic order can
be used to find solutions to systems of polynomial equations.

Recall that the ideal quotient $(I : f)$ is defined to be the ideal
$\setdef{g \in S}{f g \in I}$. Since $S$ is Noetherian, the ascending
chain of ideals $(I : f) \subset (I : f^2) \subset (I : f^3) \subset
\ldots $ is stationary; call this stationary limit $(I : f^\infty) =
\setdef{ g \in S }{ f^m g \in I \mbox{ for some } m }$.

\begin{prop}
\label{id7}
If $[x_n^{a_1} f_1, \ldots, x_n^{a_r} f_r]$ is a \std\ basis for the
ideal $I$ with respect to the reverse lexicographic order, and if none
of $f_1, \ldots, f_r$ are divisible by $x_n$, then $F = [f_1, \ldots,
f_r]$ is a \std\ basis for the ideal $(I : x_n^\infty)$. In
particular, $f_1, \ldots, f_r$ generate $(I : x_n^\infty)$.
\end{prop}

\begin{proof}
(\cite{bay82}, \cite{bs87a}) We have $F \subset (I : x_n^\infty)$. Let
$f \in (I : x_n^\infty)$; $x_n^m f \in I$ for some $m$, so $\ini(x_n^m
f)$ is a multiple of $\ini(x_n^{a_i} f_i)$ for some $i$. Since $f_i$
is not divisible by $x_n$, $\ini(f_i)$ is not divisible by $x_n$, so
$\ini(f)$ is a multiple of $\ini(f_i)$.  Thus, $F$ is a
\std\ basis for $(I : x_n^\infty)$. By \lemref{id5}, $f_1, \ldots, f_r$
generate $(I : x_n^\infty)$.
\end{proof}

If $I = \q_0 \cap \q_1 \cap \ldots \cap \q_t$ is a primary
decomposition of $I$, then $(I : x_n^\infty) = ( \cap \q_i :
x_n^\infty) = \cap (\q_i : x_n^\infty)$. We have $(\q_i : x_n^\infty)
= (1)$ if the associated prime $\p_i$ of $\q_i$ contains $x_n$, and
$(\q_i : x_n^\infty) = \q_i$ otherwise.  Thus, if $I$ defines the
subscheme $X \subset \PP^n$, then $(I : x_n^\infty)$ defines the
subscheme consisting of those primary components of $X$ not supported
on the hyperplane $\{x_n = 0\}$.

$(I : x_n^\infty)$ is saturated, because it cannot have $(x_0, \ldots,
x_n)$ as an associated prime. If $x_n$ belongs to none of the
associated primes of $I$ except $(x_0, \ldots, x_n)$, or equivalently
if $\{x_n = 0\}$ is a generic hyperplane section of $X \subset \PP^n$,
then $(I : x_n^\infty) = \sat I$. Thus, the reverse lexicographic
order can be used to find the saturation of $I$.

One of the most important uses of \std\ bases is that they lead to
canonical representations of polynomials modulo an ideal $I$, i.e.  a
division algorithm in which every $f \in S$ is written canonically as
$f = \sum g_i f_i + h$, where $[f_1, \ldots, f_r]$ is a \std\ basis
for $I$, and $h$ is the remainder after division.

Recall the division algorithm for inhomogeneous, univariate
polynomials $f(x)$, $g(x) \in k[x]$: Let $\ini(f)$ denote the highest
degree term of $f$. The remainder of $g$ under division by $f$ can be
recursively defined by $$R_f(g) = R_f(g -c x^a f)$$ if $\ini(f)$
divides $\ini(g)$, where $c x^a = \ini(g) / \ini(f)$, and by $$R_f(g)
= g$$ otherwise.

Division can be generalized to homogeneous polynomials $f_1, \ldots,
f_r, g \in S$, given a multiplicative order on $S$ (\cite{hir64},
\cite{bri73}, \cite{gal74}, \cite{sch80}): The remainder $R_F (g)$ of $g$
under division by the list of polynomials $F = [f_1, \ldots, f_r]$ can
be recursively defined by $$R_F(g) = R_F(g - c \x^A f_i)$$ for the
least $i$ so $\ini(g)$ is a multiple $c \x^A$ of $\ini(f_i)$, and by
$$R_F(g) = \ini(g) + R_F(g - \ini(g))$$ if $\ini(g)$ is not a multiple
of any $\ini(f_i)$. $R_F(g)$ is an element of $S$.

Thus, the fate of $\ini(g)$ depends on whether or not $\ini(g) \in
(\ini(f _1), \ldots, \ini(f_r))$. Let $I$ be the ideal generated by
$f_1, \ldots, f_r$.  If $F = [f_1, \ldots, f_r]$ fails to be a \std\
basis for $I$, then the remainder is poorly behaved. For example, with
respect to the lexicographic order on $k[x, y]$, $$R_{[xy, x^2 + y^2
]}(x^2 y) = x^2 y - x(xy) = 0,$$ but $$R_{[x^2 + y^2, xy]}(x^2 y) =
x^2 y - y(x^2 + y^2) = -y^3,$$ so the remainder $R_F(g)$ is dependent
on the order of the list $F$. Note that $x^2 y \in (x^2 + y^2, xy)$.

If on the other hand, $F$ is a \std\ basis for the ideal $I$, then
$R_F(g)$ is a $k$-linear combination of monomials not belonging to
$\ini(I)$. By \corref{id3}, these monomials form a $k$-basis for $S /
I$, so each polynomial in $S$ has a unique representation in terms of
this $k$-basis, modulo the ideal $I$. The remainder gives this unique
representation, and is independent of the order of $F$ (but dependent
on the multiplicative order chosen for the monomials of $S$). In
particular, $R_F(g) = 0$ if and only if $g \in I$.

An algorithm for computing a \std\ basis for $I$ from a set of
generators for $I$ was first given by Buchberger (\cite{buc65},
\cite{buc76}). This algorithm was discovered independently by Spear
(\cite{spe77}, \cite{zac78}), Bergman \cite{ber78}, and Schreyer
\cite{sch80}. It was termed the division algorithm by Schreyer, after
the division theorem of Hironaka (\cite{hir64}, \cite{bri73},
\cite{gal74}).

Define $S(f_i, f_j)$ for $i < j$ by $$S(f_i, f_j) = b \x^B f_i - c
\x^C f_j,$$ where $\x^A = b \x^B \ini(f_i) = c \x^C \ini(f_j)$ is the
least common multiple of $\ini(f_i)$ and $\ini(f_j)$. $b \x^B f_i$ and
$c \x^C f_j$ each have $\x^A$ as lead term, so $\x^A$ cancels out in
$S(f_i, f_j)$, and $\x^A > \ini(S(f_i, f_j))$.

If $F$ is a \std\ basis for the ideal $I$, then $R_F(S(f_i, f_j)) = 0$
for each $i < j$, since $S(f_i, f_j) \in I$. Conversely,

\begin{prop}[Buchberger]
\label{id8}
If $R_F(S(f_i, f_j)) = 0$ for each $i < j$, then $F = [f_1, \ldots,
f_r]$ is a \std\ basis for the ideal $I = (f_1, \ldots, f_r)$.
\end{prop}

See \cite{buc65}, \cite{buc76}. We postpone a proof until the theory
has been extended to $S$-modules. This result can also be thought of
as an explicit converse to the assertion that if $F$ is a \std\ basis,
then division is independent of the order of $F$: Whenever we have a
choice in division between subtracting off a multiple of $f_i$ and a
multiple of $f_j$, the difference is a multiple of $S(f_i, f_j)$. If
division is independent of the order of $F$, then these differences
must have remainder zero, so by \propref{id8}, $F$ is a \std\ basis.

As sketched in \secref{geometry}, \propref{id8} can be used to compute
a \std\ basis from a set of generators $f_1, \ldots, f_r$ for the
ideal $I$: For each $i < j$ so $f_{r+1} = R_F(S(f_i, f_j)) \ne 0$,
adjoin $f_{r+1}$ to the list $F = [f_1, \ldots, f_r]$. Note that
$f_{r+1} \in I$. By iterating until no new polynomials are found, a
\std\ basis $F$ is obtained for $I$. This process terminates because
$S$ is Noetherian, and each new basis element corresponds to a
monomial not in the ideal generated by the preceding lead terms.

We now extend this theory to $S$-modules. Let $M$ be a graded,
finitely generated $S$-module, given by the exact sequence of graded
$S$-modules $$M_1 \stackrel{F}{\longrightarrow} M_0 \longrightarrow M
\longrightarrow {\bf 0},$$
where $M_0 = S e_{01} \oplus \ldots \oplus S e_{0q}$ and $M_1 = S
e_{11} \oplus \ldots \oplus S e_{1r}$ are free $S$-modules with
$\deg(e_{ij}) = d_{ij}$ for each $i$, $j$. We now think of $F$ both as
a list $[f_1, \ldots, f_r]$ of module elements, and as a map between
free modules: Let $f_i = F(e_{1i}) \ne 0$ for $i =1, \ldots, r$, and
let $I \subset M_0$ be the homogeneous submodule generated by $f_1,
\ldots, f_r$.  Thus, $M = M_0 / I$.

A monomial of $M_0$ is an element of the form $\x^A e_{0i}$; such an
element has degree $\deg(\x^A) + d_{0i}$. An order on the monomials of
$M_0$ is multiplicative if whenever $\x^A e_{0i} > \x^B e_{0j}$, then
$\x^C \x^A e_{0i} > \x^C \x^B e_{0j}$ for all $\x^C \in S$. For some
applications, such as developing a theory of \std\ bases over
quotients of $S$, one wants this order to be compatible with an order
on $S$: If $\x^A > \x^B$, then one wants $\x^A e_{0i} > \x^B e_{0i}$
for $i = 1, \ldots, r$. The orders encountered in practice invariably
satisfy this second condition, but it does not follow from the first,
and we do not require it here.

One way to extend a multiplicative order on $S$ to a compatible
multiplicative order on $M_0$ is to declare $\x^A e_{0i} > \x^B
e_{0j}$ if $i < j$, or if $i = j$ and $\x^A > \x^B$. Another way is to
assign monomials $\x^{C_1}, \ldots, \x^{C_q}$ in $S$ to the basis
elements $e_{01}, \ldots, e_{0q}$ of $M_0$, and to declare $\x^A
e_{0i} > \x^B e_{0j}$ if $\x^{A+C_i} > \x^{B+C_j}$, or if $A+C_i =
B+C_j$ and $i < j$.

Fix a choice of a multiplicative order $>$ on $M_0$. The constructions
developed for $S$ carry over intact to $M_0$, with the same proofs
(\cite{gal79}, \cite{sch80}, \cite{bay82}): Given an element $f \in
M_0$, define $\ini(f)$ to be the lead term of $f$. Define $\ini(I)$ to
be the submodule generated by the lead terms of all elements of $I
\subset M_0$; $\ini(I)$ is a monomial submodule of $M_0$ with the same
Hilbert function as $I$. Define $F = [f_1, \ldots, f_r] \subset I$ to
be a \std\ basis for $I$ if $\ini(f_1), \ldots, \ini(f_r)$ generate
$\ini(I)$; a set of generators for $I$ need not be a \std\ basis for
$I$, but a \std\ basis for $I$ generates $I$. Given an element $g \in
M_0$, define $R_F(g) \in M_0$ exactly as was done for the free module
$S$. If $F$ is a \std\ basis for $I$, then $R_F(g) = 0$ if and only if
$g \in I$.

The quotient of $g$ under division by $f_1, \ldots, f_r$ can be
recursively defined by $$Q_F(g) = c \x^A e_{1i} + Q_F(g - c \x^A
f_i)$$ for the least $i$ so $\ini(g)$ is a multiple $c \x^A$ of
$\ini(f_i)$, and by $$Q_F(g) = Q_F(g - \ini(g))$$ if $\ini(g)$ is not
a multiple of any $\ini(f_i)$. The quotient is an element of $M_1$.

Following the recursive definitions of the remainder and quotient, it
can be inductively verified that $$g = F(Q_F(g)) + R_F(g).$$ If $F$ is
a \std\ basis for $I$, and $g \in I$, then $R_F(g) = 0$, so the
quotient lifts $g$ to $M_1$. In this case, the quotient can be thought
of as expressing $g$ in terms of $f_1, \ldots, f_r$.

Define $S(f_i, f_j)$ for $i < j$ by $$S(f_i, f_j) = b \x^B f_i - c
\x^C f_j,$$ if $\ini(f_i)$ and $\ini(f_j)$ have a least common
multiple $\x^A e_{0k} = b \x^B \ini(f_i) = c \x^C \ini(f_j)$. Leave
$S(f_i, f_j)$ undefined if $\ini(f_i)$ and $\ini(f_j)$ lie in
different summands of $M_0$, and so don't have common multiples.

Recall that the module of syzygies of $f_1, \ldots, f_r$ is defined to
be the kernel of the map $F$, which is the submodule of $M_1$
consisting of all $h \in M_1$ so $F(h) = 0$. Thus, if $h = h_1 e_{11}
+ \ldots + h_r e_{1r}$ is a syzygy, then $h_1 f_1 + \ldots + h_r f_r =
0$. Let $J \subset M_1$ denote the module of syzygies of $f_1, \ldots,
f_r$, and let $K \subset M_1$ denote the module of syzygies of
$\ini(f_1), \ldots, \ini(f_r)$.

Define the map $\ini(F): M_1 \rightarrow M_0$ by $\ini(F)(e_{1i}) =
\ini(f_i)$; $K$ is the kernel of $\ini(F)$. For each $i < j$ so
$S(f_i, f_j)$ is defined, define $t_{ij}$ to be the element $$t_{ij} =
b \x^B e_{1i} - c \x^C e_{1j} \in M_1,$$ where $\x^A e_{0k} = b \x^B
\ini(f_i) = c \x^C \ini(f_j)$ is the least common multiple of
$\ini(f_i)$ and $\ini(f_j)$, as before. $\ini(F)(t_{ij}) = 0$, so each
$t_{ij}$ belongs to the syzygy module $K$. Observe that $F(t_{ij}) =
S(f_i, f_j)$.

Assign the following multiplicative order on $M_1$, starting from the
order on $M_0$ (\cite{sch80}; see also \cite{mm86}): Let $\x^A e_{1i}
> \x^B e_{1j}$ if $\x^A \ini(f_i) > \x^B \ini(f_j)$, or if these terms
are $k$-multiples of each other and $i < j$. If the order on $M_0$ is
compatible with an order on $S$, then this order on $M_1$ is
compatible with the same order on $S$.

With respect to this order on $M_1$, we have

\begin{lemma}
\label{tlem}
The list $[\,t_{ij}\,]$ is a \std\ basis for the module $K$ of
syzygies of $\ini(f_1), \ldots, \ini(f_r)$.
\end{lemma}

\begin{proof}
Let $h \in M_1$, so $\ini(F)(h) = 0$. Then $\ini(F)(\ini(h))$ is
canceled by $\ini(F)(h - \ini(h))$ in $M_0$. Therefore, if $\ini(h) =
\x^A e_{1i}$, then $h$ has another term $\x^B e_{1j}$ so $\x^A
\ini(f_i)$ and $\x^B \ini(f_j)$ are $k$-multiples of each other and $i
< j$. Thus, $t_{ij}$ is defined and $\ini(t_{ij})$ divides $\ini(h)$,
so $[\,t_{ij}\,]$ is a \std\ basis for $K$.
\end{proof}

Thus, the set $\{t_{ij}\}$ generates $K$. In general, the
$[\,t_{ij}\,]$ are far from being a minimal \std\ basis for $K$; we
consider the effects of trimming this list in \propref{id9} below.

Define $$s_{ij} = t_{ij} - Q_F(S(f_i, f_j))$$ whenever $R_F(S(f_i,
f_j)) = 0$. Note that $\ini(s_{ij}) = \ini(t_{ij})$. Each $s_{ij}$ is
the difference of two distinct elements of $M_1$, each of which is
mapped by $F$ to $S(f_i, f_j)$, so $F(s_{ij}) = 0$. In other words,
$s_{ij}$ belongs to the syzygy module $J$. Conversely,

\begin{prop}[Richman, Spear, Schreyer]
\label{id9}
Choose a set of pairs $T = \{(i, j)\}$ such that the set
$\{t_{ij}\}_{(i, j) \in T}$ generates the module $K$ of syzygies of
$\ini(f_1), \ldots, \ini(f_r)$. If $R_F(S(f_i, f_j)) = 0$ for each
$(i, j)\in T$, then

(a) $F = [f_1, \ldots, f_r]$ is a \std\ basis for $I$;

(b) the set $\{s_{ij}\}_{(i, j) \in T}$ generates the module $J$ of
syzygies of $f_1, \ldots, f_r$.

\noindent Moreover,

(c) if $[\,t_{ij}\,]_{(i,j)\in T}$ is a \std\ basis for $K$, then
$[\,s_{ij}\,]_{(i,j)\in T}$ is a \std\ basis for $J$.
\end{prop}

\begin{proof}
(\cite{ric74}, \cite{spe77}, \cite{zac78}, \cite{sch80}) First,
suppose that $[\,t_{ij}\,]_{(i,j) \in T}$ is a \std\ basis for $K$.
Let $h \in J$, so $F(h) = 0$. By the same reasoning as in the proof of
\lemref{tlem}, we can find $(i,j) \in T$ so $\ini(t_{ij})$ divides
$\ini(h)$.  Since $\ini(s_{ij}) = \ini(t_{ij})$, $\ini(s_{ij})$ also
divides $\ini(h)$, so $[\,s_{ij}\,]_{(i,j) \in T}$ is a \std\ basis
for $J$, proving (c).

Now, suppose that $\{t_{ij}\}_{(i,j) \in T}$ merely generates $K$. Let
$T^\prime$ be a set of pairs so $[\,t_{\ell m}\,]_{(\ell,m) \in
T^\prime}$ is a \std\ basis for $K$. It is enough to construct a list
$[\,u_{\ell m}\,]_{(\ell,m) \in T^\prime}$ of elements of $J$,
generated by $\{s_{ij}\}_{(i,j) \in T}$, so $\ini(u_{\ell m}) =
\ini(t_{\ell m})$ for all $(\ell,m) \in T^\prime$. Then by the
preceding argument, $[\,u_{\ell m}\,]_{(\ell,m) \in T^\prime}$ is a
\std\ basis for $J$, so $\{s_{ij}\}_{(i,j) \in T}$ generates $J$.

Write each $t_{\ell m} = \sum g_{\ell m i j} t_{ij}$, for $(\ell,m)
\in T^\prime$ and $(i,j) \in T$, in such a way that the terms of
$t_{\ell m}$ and each term of each product $g_{\ell m i j} t_{ij}$ map
via $\ini(F)$ to multiples of the same monomial in $M_0$. In other
words, find a minimal expression for each $t_{\ell m}$, which avoids
unnecessary cancellation.  Then define $$u_{\ell m} = \sum g_{\ell m i
j} s_{ij}.$$ We have $\ini(u_{\ell m}) = \ini(t_{\ell m})$, proving
(b).

Let $f \in I$, and choose $g \in M_1$ so $f = F(g)$. Let $h \in M_1$
be the remainder of $g$ under division by $[\,u_{\ell m}\,]_{(\ell,m)
\in T^\prime}$; $f = F(h)$. Since $\ini(h)$ is not a multiple of any
$\ini(u_{\ell m}) = \ini(t_{\ell m})$, the lead term of $F(\ini(h))$
is not canceled by any term of $F(h - \ini(h))$. Therefore, if
$\ini(h) = a \x^A e_{1i}$, then $\ini(f_i)$ divides $\ini(F)$. Thus,
$F = [f_1, \ldots, f_r]$ is a \std\ basis for $I$, proving (a).
\end{proof}

\propref{id8} follows as a special case of this result.

The above proof can be understood in terms of an intermediate initial
form $\ini_0(h)$ for $h \in M_1$: Apply the map $\ini(F)$ separately
to each term of $h$, and let $\x^A \in M_0$ be the greatest monomial
that occurs in the set of image terms. Define $\ini_0(h)$ to be the
sum of all terms of $h$ which map via $\ini(F)$ to multiples of
$\x^A$. Then $\ini$ refines $\ini_0$, for according to the order we
have defined on $M_1$, $\ini(h)$ is the term of $\ini_0(h)$ lying in
the summand of $M_1$ whose basis element $e_i$ has the smallest index
$i$.

In this language, $t_{ij} = \ini_0(t_{ij}) = \ini_0(s_{ij})$. Our
expressions for the $t_{\ell m}$ have the property that each $g_{\ell
m i j} t_{ij} = \ini_0(g_{\ell m i j} t_{ij})$, with each term of each
product for a given $t_{\ell m}$ mapping via $\ini(F)$ to multiples of
the same monomial $\x^A$. Thus, each $\ini_0(g_{\ell m i j} s_{ij}) =
g_{\ell m i j} t_{ij}$; the tails $g_{\ell m i j} (s_{ij} -
\ini_0(s_{ij}))$ stay out of our way, mapping termwise via $\ini(F)$
to monomials which are less than $\x^A$ with respect to the order on
$M_0$.

Observe that $Q_F(g)$ is a linear combination of monomials not
belonging to $\ini(J)$, for any $g \in M_0$.

In \cite{buc79}, Buchberger gives a criterion for selecting a set $T$
of pairs $(i, j)$ in the case where $I$ is an ideal: If $(i_0, i_1),
(i_1, i_2), \ldots, (i_{s-1}, i_s) \in T$, and the least common
multiple of $\ini(f_{i_0}), \ini(f_{i_1}), \ldots, \ini(f_{i_s})$ is
equal to the least common multiple of $\ini(f_{i_0})$ and
$\ini(f_{i_s})$, then $(i_0, i_s)$ need not belong to $T$.  In other
words, if $t_{i_0 i_s} \in (t_{i_0 i_1}, \ldots, t_{i_{s-1} i_s})$,
then the pair $(i_0, i_s)$ is unnecessary; this condition is
equivalent to the condition of \propref{id9}, for the case of an
ideal.

Suppose that we wish to compute the syzygies of a given set of
elements $g_1, \ldots, g_s$ of $M_0$. To do this, compute a \std\
basis $f_1, \ldots, f_r$ for the submodule $I \subset M_0$ generated
by $g_1, \ldots, g_s$.  Keep track of how to write each $f_i$ in terms
of $g_1, \ldots, g_s$. Using these expressions, each syzygy of $f_1,
\ldots, f_r$ can be mapped to a syzygy of $g_1, \ldots, g_s$. These
images generate the module of syzygies of $g_1, \ldots, g_s$; the set
of syzygies obtained in this way is not in general minimal.

Syzygies can be used to find a minimal set of generators for a
submodule $I \subset M_0$ from a given set of generators $g_1, \ldots,
g_s$: If $h_1 g_1 + \ldots + h_r g_r = 0$ is a syzygy of $g_1, \ldots,
g_s$ with $h_1 \in k$, then $g_1 = (h_2 g_2 + \ldots + h_r g_r) /
h_1$, so $g_1$ is not needed to generate $I$. All unnecessary
generators can be removed in this way.

Alternatively, a careful implementation of \std\ bases can directly
find minimal sets of generators for submodules: Starting from an
arbitrary set of generators, we can eliminate unnecessary
generators degree by degree, by removing those which reduce to
zero under division by a \std\ basis for the ideal generated by
the preceding generators.

Either way, we can trim the set of syzygies computed via \std\ bases
for a given set of generators $g_1, \ldots, g_s$ of $I$, to obtain a
minimal set of generators for the syzygy module $J$. By starting with
a minimal generating set for $I$, and iterating this method, a minimal
free resolution can be found for $I$.

A beautiful application of these ideas yields a proof of the Hilbert
syzygy theorem, that minimal free resolutions terminate
(Schreyer \cite{sch80}, \cite{sch91}, for an exposition see also
Eisenbud \cite{eis92}). At each stage of a resolution, order the \std\
basis $F$ for $I$ in such a way that for each $i < j$, letting $\ini(f_i) = a
\x^A e_{0k}$ and $\ini(f_j) = b \x^B e_{0\ell}$, we have $\x^A > \x^B$
in the lexicographic order. If the variables $x_1,\ldots, x_m$ are
missing from the initial terms of the $f_i$, then the variables
$x_1,\ldots, x_{m+1}$ will be missing from the initial terms of the
syzygies $s_{ij}$. Iterating, we run out of variables, so the resolution
terminates.

\section{Bounds}
\label{bounds}
How hard are the algorithms in algebraic geometry? We describe some key bounds.
The best known example is the bound established by G. Hermann \cite{her26} for
ideal membership:

\begin{thm}[G. Hermann]\label{thm1}
Let $k$ be any field, let $(f_1,...,f_k) \subset k[x_1,...,x_n]$ and
let $d =
\max(\deg(f_i))$. If $g \in (f_1,...,f_k)$, then there is an expression
$$
g = \sum_{i=1}^k a_i f_i
$$
where $\deg(a_i) \le \deg(g) + 2(kd)^{2^{n-1}}$.
\end{thm}

This type of bound is called \quotes{doubly exponential}. However,
with the advent of the concept of coherent sheaf cohomology
\cite{ser55} and the systematic study of vanishing theorems, it has
become apparent that the vanishing of these groups in high degrees is
almost always the most fundamental bound. The concept of an ideal
being \quotes{m-regular} or \quotes{regular in degrees $\ge m$} was
introduced by one of us \cite{mum66} by generalizing ideas of
Castelnuovo:

\begin{defn}\label{def1}
\footnote{
The definition has been slightly modified so as to apply to ideals $I$
instead of the corresponding sheaf of ideals $\Ih$.  } Let $k$ be any
field, let $I \subset k[x_0,\ldots,x_n]$ be an ideal generated by
homogeneous polyomials, let $I_d$ be the homogeneous elements in $I$
of degree $d$, let $I$ be the corresponding sheaf of ideals in
$\Oh_{\proj n}$, and let $I(d)$ be the $d^{\it th}$ twist of $I$.
Then the following properties are equivalent and define the term
\quotes{m-regular}:

(a) the natural map $I_m \rightarrow H^0(\Ih(m))$ is an isomorphism
and \mbox{$H^i(\Ih(m-i))$} $ = (0)$, $1 \le i \le n$

(b) the natural maps $I_d \rightarrow H^0(\Ih(d))$ are isomorphisms
for all $d \ge m$ and $H^i(\Ih(d))$ = $(0)$ if $d+i \ge m$, $i \ge 1$.

(c)  Take a minimal resolution of $I$ by free graded $k[X]$-modules:
$$
0 \, \rightarrow \,
\Limits{\oplus}{\alpha=1}{r_n} \,
k[\x] \cdot e_{\alpha,n} \,
\stackrel{\phi_n}{\rightarrow} \,
\ldots
\stackrel{\phi_1}{\rightarrow} \,
\Limits{\oplus}{\alpha=1}{r_0} \,
k[\x] \cdot e_{\alpha,0} \,
\stackrel{\phi_0}{\rightarrow} \,
k[\x] \longrightarrow k[\x]/I \rightarrow 0.  $$ Then
$\deg(e_{\alpha,i}) \le m+i$ for all $\alpha$, $i$. (In particular, if
$f_\alpha = \phi_0(e_{\alpha,0})$, then $f_1,\ldots,f_{r_0}$ are
minimal generators of $I$, and $\deg(e_{\alpha,0}) = \deg(f_\alpha)
\le m$.)
\end{defn}

The intuitive idea is that past degree $m$, nothing tricky happens in
the ideal $I$.  Unfortunately, neither (a), (b) nor (c) can be
verified by any obvious finite algorithm.  This lack of a finitely
verifiable criterion for $m$-regularity has been remedied by a joint
result of the first author and M. Stillman \cite{bs87a}:

\begin{thm}[Bayer-Stillman]\label{thm2}
$I$ is $m$-regular if and only if the degrees of the minimal set of
generators of $I$ are at most $m$, and there exists a set
$y_0,\ldots,y_\ell$ of linear combinations of $x_0,...,x_n$ such that
for all homogeneous $f$ of degree $m$,
\begin{eqnarray*}
y_0 f \in I &\Rightarrow& f \in I \\ y_1 f \in I &\Rightarrow& f \in I
+ k[\x] \cdot y_0 \\ &\cdots& \\ y_\ell f \in I &\Rightarrow& f \in I
+ \sum_{i=0}^{\ell-1} k[\x] \cdot y_i
\end{eqnarray*}
and $$f \in I + \sum_{i=0}^{\ell} k[\x] \cdot y_i.$$ Moreover, if this
holds at all, it holds for $y_0,\ldots,y_\ell$ taken arbitrarily from
a Zariski-open set in the space of $\ell+1$ linear forms.
\end{thm}

To see why $m$-regularity is a key bound, we want to show that it
controls some of the geometric features of the ideal $I$.  Let's
introduce several refined notions of the \quotes{degree} of $I$:

\begin{defn}\label{def2}
If $I = \q_0 \,\cap\, \q_1 \,\cap\, \ldots \,\cap\, \q_t$ is a primary
decomposition of $I$, $\sqrt{\q_i} = \p_i$ is prime and $V(\p_i)$ is
the subvariety $Z_i$ of $\proj n$ for $i \ge 1$, while $\p_0 =
(x_0,\ldots,x_n)$ (so that $V(\p_0)= \emptyset$), then first let
$\q_1,\ldots,\q_s$ be the isolated components, (i.e., $Z_i \not\subset
Z_j$ if $1\le i\le s$, $1\le j \le t$, $i \ne j$, or equivalently,
$V(I)=Z_1 \,\cup\, \ldots \,\cup\, Z_s$ is set-theoretically the
minimal decomposition of $V(I)$ into varieties).  Then let
\begin{eqnarray*} \mult(\q_i) &=& \mbox{length $\ell$ of a maximal
chain of $\p_i$-primary ideals:} \\ & &\q_i = J_\ell \ssubset
J_{\ell-1} \ssubset \ldots
\ssubset J_1 = \p_i
\end{eqnarray*}
(Equivalently, this is the length of the local ring $k[\x]_{\p_i}/I
k[\x]_{\p_i}$, or, in the language of schemes, if $\eta$ is the
generic point of $Z_i$, then this is the length of $\Oh_{\eta,\proj
n}.$)
\begin{eqnarray*}
\deg (Z_i) &=& \mbox{usual geometric degree of $Z_i$:} \\
& &\mbox{the cardinality of $Z_i \cap L$ for almost all} \\
& & \mbox{linear spaces $L$ of complementary dimension.} \\
\geomdeg_r(I) &=&
\sum_{\stackrel{\mbox{\it \scriptsize $i$ such that $\dim Z_i = r$}}
{\mbox{\scriptsize $1 \le i \le s$}}} \mult(\q_i) \deg(Z_i)
\end{eqnarray*}
\end{defn}

If $\q_i$ is one of the non-isolated, or embedded components, then we
extend the concept of multiplicity more carefully: Let $$ I_i =
\left\{ \cap \q_j \mbox{\Large $\mid$} j \mbox{ such that } \p_j
\ssubset \p_i \mbox{ or equivalently } Z_j \ssupset Z_i \right\} \cap \p_i
$$ and
\begin{eqnarray*}
\mult_I(\q_i) &=& \mbox{length $\ell$ of a maximal chain of ideals:} \\
& & \q_i \cap I_i = J_\ell \ssubset J_{\ell-1} \ssubset \ldots
\ssubset J_0 = I_i \\ & & \mbox{where each $J_k$ satisfies: } ab \in
J_k, a \not\in \p_i \Rightarrow b \in J_k.
\end{eqnarray*}
(Equivalently, $J_k$ equals $\q_k \cap I_i$ for some $\p_i$-primary
ideal $\q_k$.)
In particular:
\begin{eqnarray*}
I_0 &=& \Bigcap{j=1}{t} \q_j \mbox{ is known as } \sat I \mbox{, and}
\\
\mult_I(\q_0) &=& \mbox{length $\ell$ of a maximal chain of ideals} \\
& & I = J_\ell \ssubset J_{\ell-1} \ssubset \ldots \ssubset J_0 = \sat
I \\ &=& \dim_k(\sat I/I).
\end{eqnarray*}
For $s+1 \le i \le t$, an equivalent way to define $\mult_I(\q_i)$ is
as the length of the module $$I_i k[\x]_{\p_i}/I k[\x]_{\p_i}$$ or, in
the language of schemes, the length of $$I_i{\Oh}_{\eta,\proj
n}/I{\Oh}_{\eta,\proj n}$$ where $\eta$ is the generic point of $Z_i$.

Then write
$$\arithdeg_r(I) =
\sum_{\stackrel{\mbox{\it \scriptsize $i$ such that $\dim Z_i = r$}}
{\mbox{\scriptsize $1 \le i \le s$}}}
\mult_I(\q_i) \deg(Z_i)$$
and
$$\arithdeg_{-1}(I) = \mult_I(\q_0).$$

The idea here is best illustrated by an example:  let
$$I = (x_1^2,x_1 x_2) \subset k[x_0,x_1,x_2].$$
Then
\begin{eqnarray*}
I &=& \q_1 \cap \q_2 \\ q_1 &=& (x_1),\; \p_1 = (x_1),\; Z_1 =
\{\mbox{line } x_1=0\} \\ q_2 &=& (x_1^2,x_1x_2,x_2^2),\; \p_2 =
(x_1,x_2), \; Z_2 = \{\mbox{point } (1,0,0)\}.
\end{eqnarray*}
Then
$$\deg(Z_1)=1,\; \mult(q_1) = 1$$
so
$$\geomdeg_1(I) = \arithdeg_1(I) = 1.$$

One might be tempted to simply define $$\mult_I(\q_2) = \mbox{ length
of chain of $\p_2$-primary ideals between } \q_2, \p_2$$ and since
$$k[\x]_{\q_2}/\q_2 k[\x]_{\p_2} \cong K \cdot 1+K \cdot x_1 + K \cdot
x_2, K=k(x_0)$$ this is $3$.  But embedded components are not unique!
In fact,
\begin{eqnarray*}
I &=& \q_1 \cap \q_2^\prime \\
\q_2^\prime &=& (x_1^2 x_2) \mbox{ also},
\end{eqnarray*}
which leads to $$k[\x]_{\p_2}/\q_2^\prime k[\x]_{\p_2} \cong K \cdot
1+ K \cdot x_2$$ which has length $2$.  The canonical object is not
the local ring $k[\x]_{\p_2}/\q_2 k[\x]{\p_2}$ but the ideal
$$\mbox{\rm Ker} \left( k[\x]_{\p_2} / I k[\x]_{\p_2}
\rightarrow
k[\x]_{\p_2} / \p_2 k[\x]_{\p_2} \right)
\cong k \cdot x_1$$
which has length $1$.  Thus, the correct numbers are
$$\mult_I(\q_2) = 1$$
and
\begin{eqnarray*}
\geomdeg_0(I) &=& 0 \\
\arithdeg_0(I) &=&1.
\end{eqnarray*}

Now the question arises: find bounds on these degrees in terms of
generators of $I$.  For geometric degrees, a straightforward extension
of Bezout's theorem gives:

\begin{prop}\label{prop3}
Let $d(I)$ be the maximum of the degrees of a minimal set of
generators of $I$.  Then $$\geomdeg_r(I) \le d(I)^{n-r}.$$
\end{prop}

A proof can be found in \cite{mw83}.  The idea is clear from a simple
case: Suppose $f, g, h \in K[x,y,z]$ and $f=g=h=0$ consists of a curve
$C$ and $\ell$ points $P_i$ off $C$.  We can bound $\ell$ like this:
Choose $2$ generic combinations $f^\prime$, $g^\prime$ of $f$, $g$,
$h$ so that $f^\prime=g^\prime=0$ does not contain a surface.  It must
be of the form $C \cup C^\prime$, $C'$ one-dimensional, containing all
the $P_i$ but not the generic point of $C$.  Then by the usual Bezout
theorem $$\deg C^\prime \le
\deg f^\prime \;\deg g^\prime=d(I)^2.$$ Let $h^\prime$ be a
$3$\thuh{rd} generic combination of $f,g,h$.  Then $C^\prime \cap
\{h^\prime = 0\}$ consists of a finite set of points including the
$P_i$'s.  Thus
\begin{eqnarray*}
\ell &=& \# P_i \le \#(C^\prime \cap \{h^\prime = 0\}) \\
& \le & \deg C^\prime \cdot d(I) \mbox{ by Bezout's theorem} \\
& \le & d(I)^3.
\end{eqnarray*}
Can $\arithdeg(I)$ be bounded in the same way?  In fact, it cannot, as
we will show below.  Instead, we have

\begin{prop}\label{prop4}
If $m(I)$ is the regularity of $I$, then for $-1 \le r \le n$,
$$\arithdeg_r(I) \le
{m(I)+n-r-1 \choose n-r}
\le m(I)^{n-r}$$
\end{prop}
which replaces $d(I)$ by the regularity of $I$.  A proof is given in
the technical appendix.

We have introduced two measures of the complexity of a homogeneous
ideal $I$.  The first is $d(I)$, the maximum degree of a polynomial in
a minimum set of generators of $I$.  The second is $m(I)$, which
bounds the degrees of generators and of all higher order syzygies in
the resolution of $I$ (\defref{def1} (c)).  Obviously, $$d(I) \le
m(I).$$

A very important question is how much bigger can $m(I)$ be than
$d(I)$?  The nature of the answer was conjectured by one of us in his
thesis \cite{bay82} and this conjecture is being borne out by
subsequent investigations.  This conjecture is that in the worst case
$m(I)$ is roughly the $(2^n)$\thuh{th} power of $d(I)$ -- a bound like
G. Hermann's. But that if $I = I(Z)$ where $Z$ is geometrically nice,
e.g. is a smooth irreducible variety, then $m(I)$ is much smaller,
like the $n$\thuh{th} power of $d(I)$ or better. This conjecture then
has three aspects:

\begin{tabular}{@{\hspace{5pt}}ll}
(1) & a doubly exponential bound for $m(I)$ in terms of $d(I)$, \\ &
which is always valid,\\ (2) & examples of $I$ where the bound in (1)
is best possible, or nearly so,\\ (3) &much better bounds for $m(I)$
\\ & valid if $V(I)$ satisfies various conditions.\\
\end{tabular}

All three aspects are partially proven, but none are completely
clarified yet. We will take them up one at a time.

A doubly exponential bound for $m(I)$ in terms of $d(I)$ may be
deduced easily {\em in characteristic zero} from the work of M. Giusti
\cite{giu84} and A. Galligo
\cite{gal79}:

\begin{thm}\label{thm4}
If $\mbox{\rm char}(k) = 0$ and $I \subset k[x_0, \ldots x_n]$ is any
homogeneous ideal, then $$m(I) \le (2 d(I))^{2^{n-1}}.$$
\end{thm}

It seems likely that \thmref{thm4} holds in characteristic $p$,
too. A weaker result can be derived quickly in any characteristic by
straightforward cohomological methods:

\begin{prop}\label{prop5}
If $I \subset k[x_0, \ldots x_n]$ is any homogeneous ideal, then
$$m(I) \le (2 d(I))^{n!}.$$
\end{prop}

The proof is given in the technical appendix.

Next, we ask whether \thmref{thm4} is the best possible, or
nearly so. The answer is yes, because of a very remarkable example due
to E. Mayr and A. Meyer
\cite{mm82}.

\begin{example}
Let $I_n^A$ be the ideal in $10 n$ variables $S^{(m)}, F^{(m)},
C_i^{(m)}, B_i^{(m)}, 1 \le i \le 4, 1 \le m \le n$ defined by the $10
n - 6$ generators $$\begin{array}{rl} 2 \le m \le n \;\; &\left\{\;\;
\begin{array}{l}
S^{(m)} - S^{(m-1)} C_1^{(m-1)} \\ F^{(m)} - S^{(m-1)} C_4^{(m-1)} \\
C_i^{(m)}F^{(m-1)}B_2^{(m-1)} -
C_i^{(m)}B_i^{(m)}F^{(m-1)}B_3^{(m-1)}, \;\;1 \le i
\le 4 \end{array}
\right.\\
\\
1 \le m \le n-1 \;\; &\left\{\;\;
\begin{array}{l}
F^{(m)}C_1^{(m)}B_1^{(m)} - S^{(m)} C_2^{(m)} \\
F^{(m)}C_2^{(m)} - F^{(m)} C_3^{(m)} \\
S^{(m)}C_3^{(m)}B_1^{(m)} - S^{(m)}C_2^{(m)}B_4^{(m)} \\
S^{(m)}C_3^{(m)} - F^{(m)}C_4^{(m)}B_4^{(m)}
\end{array}
\right.\\
\\
&\mbox{\hskip 20pt}C_i^{(1)}S^{(1)} - C_i^{(1)}F^{(1)}(B_i^{(1)})^2,
\;\;1 \le i \le 4
\end{array}$$
Let $I_n^H$ be the ideal gotten from $I_n^A$ by homogenizing with an
extra variable $u$. Then Mayr and Meyer \cite[lemma 8, p. 318]{mm82}
prove:
\end{example}

\begin{lemma}
 Let $e_n=2^{2^n}$. If $M$ is any monomial in these variables,
$S^{(n)}C_i^{(n)}-F^{(n)} M \in I_n^A$ if and only if
$$M=C_i^{(n)}(B_i^{(n)})^{e_n},$$
and $S^{(n)}C_i^{(n)} - S^{(n)} M \in I_n^A$ if and only if
$$M=C_i^{(n)}.$$
\end{lemma}

Now note that the generators of $I_n^A$ and $I_n^H$ are all of the
very simple type given by a difference of two monomials.  Quite
generally, if
\begin{eqnarray*}
J &\subset& k[x_1,\ldots,x_n] \\
J &=& (\ldots, \x^{\alpha_i} - \x^{\beta_i}, \ldots )_{1 \le i \le k}
\end{eqnarray*}
then the quotient ring $k[\x]/J$ has a very simple form.  In fact, we
get an equivalence relation between monomials generated by
$$\x^{\alpha_i+\gamma} \sim \x^{\beta_i+\gamma}, \mbox{ any } i,
\gamma$$ and $$k[\x]/J \cong \oplus_\delta \; k \cdot \x^\delta$$
where $\delta$ runs over a set of representatives of each equivalence
class.

Bearing this in mind, let's look at the $1$\thuh{st} order syzygies
for the homogeneous ideal: $$J_n^H = (S^{(n)},F^{(n)},I_n^H).$$
$S^{(n)}$ and $F^{(n)}$ are part of a minimal set of generators, and
let $f_\alpha \in I_n^H$ complete them.  Then syzygies are equations:
$$p \, S^{(n)} \, + \, q \, F^{(n)} \, + \, \sum \, r_\alpha \,
f_\alpha \, = \, 0.$$ One such is given by: $$\left[ u^{e_n+e} \,
C_i^{(n)} \right] \, S^{(n)} \, + \,
\left[ - u^e \, (B_i^{(n)})^{e_n} \, C_i^{(n)} \right] \, F^{(n)} \,
+ \, \sum \, R_\alpha \, f_\alpha \, = \, 0$$ for some $R_\alpha$, and
some $e \ge 0$ (the extra power $u^e$ is necessary because some terms
$R_\alpha f_\alpha$ have degree greater than $e_n+2$) whose degree is
$2+e_n+e$.  Now express this syzygy as a combination of a minimal set
of syzygies.  This gives us in particular:
\begin{eqnarray*}
u^{e_n+e} \, C_i^{(n)} \, &=& \, \sum \, a_\lambda \, p_\lambda \\ -
u^e \, (B_i^{(n)})^{e_n} \, C_i^{(n)} \, &=& \, \sum \, a_\lambda \,
q_\lambda \\ p_\lambda \, S^{(n)} \, + \, q_\lambda \, F^{(n)} \, + \,
\sum \, R_{\alpha\lambda} \, f_\alpha \, &=& \, 0 \, .
\end{eqnarray*}

Then for some $\lambda$, $p_\lambda$ must have a term of the form
$u^\ell$ or $u^\ell \, C_i^{(n)}$, hence the monomial $u_\ell \,
S^{(n)}$ or $u_\ell \, C_i^{(n)} \, S^{(n)}$ occurs in $p_\lambda \,
S^{(n)}$.  But by the general remark on quotient rings by such simple
ideals, this means that this term must equal some second term $M \,
S^{(n)}$ ($M$ a monomial in $p_\lambda$) or $M \, F^{(n)}$ ($M$ a
monomial in $q_\lambda$) mod $I_n^H$.  By the lemma, the first doesn't
happen and the second only happens if the term $u^\ell \, C_i^{(n)} \,
(B_i^{(n)})^{e_n}$ occurs in $q_\lambda$, in which case $e_n+1\le \deg
\, q_\lambda \,= \,
\deg(\mbox{syzygy}(p_\lambda,q_\lambda,R_{\alpha\lambda})) - 1$.  This
proves:

\begin{prop}\label{prop7}
$J_n^H$ has for its bounds:
\begin{eqnarray*}
d(J) &=& 4 \\
m(J) &\ge& 2^{2^n}+1.
\end{eqnarray*}
\end{prop}
Going on to the $3$\thuh{rd} aspect of the conjecture, consider
results giving better bounds for $m(I)$ under restrictive hypotheses
on $V(I)$.

\begin{thm}\label{thm8}
If $Z \subset \proj n$ is a reduced subscheme purely of dimension $r$,
and $I=I(Z)$ is the full ideal of functions vanishing on $Z$, then

(a) if $r \leq 1$, or $Z$ is smooth, $\hbox{char}(k)=0$ and $r \leq
3$, then: $$m(I) \le \deg Z - n + r + 1$$

(b) if char$(k)=0$ and $Z$ is smooth,
$$m(I) \le (r+1)(\deg(Z)-2)+2.$$
\end{thm}

Since $\deg(Z) \le d(I)^{n-r}$ (\propref{prop3}), these bound $m(I)$
in terms of $d(I).$

Part (a) of this are due to Gruson-Lazarsfeld-Peskine \cite{glp83} for
$r \leq 1$, and to Pinkham \cite{pin86}, Lazarsfeld \cite{laz87}, and
Ran \cite{ran90} for $r \leq 3$. It is {\it conjectured} by Eisenbud and Goto
\cite{eg84}, and others, that the bound in (a) holds for all reduced
irreducible $Z$,
and it might well hold even for reduced equidimensional $Z$ which are connected
in
codimension 1.  As this problem is now understood, the needed
cohomological arguments follow formally, once one can control the singularities
of a
projection of the variety. These singularities become progressively harder to
subdue
as the dimension of the variety increases, and are what impedes definitive
progress
beyond dimension 3.

Part (b) is due to the second author and is proven in the technical
appendix.  It has been generalized by Bertram, Ein, and Lazarsfeld
\cite{bel91} to show that any smooth characteristic $0$ variety of
codimension $e$ defined as a subscheme of ${\bf P}^n$ by hypersurfaces
of degrees $d_1 \ge \ldots \ge d_m$ is $(d_1 + \ldots d_e - e
+1)$-regular.  Since we cannot decide the previous conjecture, this is a
result of considerable practical importance, for it strongly bounds
the complexity of computing \std\ bases of smooth characterisitic $0$
varieties in terms of the degrees of the input equations.

The biggest missing link in this story is a decent bound on $m(I)$ for
any reduced equidimensional ideal $I$. We would conjecture that if a
linear bound as in part (a) doesn't hold, at the least a so-called
``single exponential'' bound, i.e. $m(I) \leq d^{0(n)}$ ought to hold.
This is an essential ingredient in analyzing the worst-case behavior
of all algorithms based on \std\ bases, and would complete the story
about what causes the bad examples discussed above.  At least in some
cases Ravi \cite{rav90} has proven that the regularity of the radical
of a scheme is no greater than the regularity of the scheme itself.

There is a direct link between the bounds that we have given so far
and the G. Hermann bound with which we started the section.  This
results from the following:

\begin{prop}\label{prop9}
Let $I^A \subset k[x_1,\ldots,x_n]$ have generators $f_1,\ldots,f_k$
and let $I^H \subset k[x_0,x_1,\ldots,x_n]$ be the ideal generated by
homogenizations $f_1^h,\ldots,f_k^h$ of the $f_i$.  Let $I^H = \q_0
\cap \ldots \cap \q_t$ be the primary decomposition of $I^H$, let
$Z_i=V(\q_i)$ and let $$\mult_\infty(I^H) = \max \left[
\mult_I(\q_{i_1}) +\ldots+\mult_I(\q_{i_k})+\mult_I(\q_0) \right]$$
where the $\max$ is taken over chains
$V((x_0)) \supset Z_{i_1} \ssupset \ldots \ssupset Z_{i_k}$.
If $g \in I^A$, then we can write:
$$g = \sum_{i=1}^k \, a_i \, f_i$$
where
$$\deg a_i \le \deg g + \mult_\infty(I^H).$$
\end{prop}

The proof goes like this: Let $g^h$ be the homogenization of $g$.
Consider the least integer $m$ such that $x_0^m g^h \in I^H$.  Since
$g \in I$, this $m$ is finite. Moreover, if $$x_0^{m} g^h = \sum
x_0^{m_i} a_i^h f_i^h$$ then $$g= \sum a_i f_i$$ and $$\deg a_i =
\deg(a^h) \le \deg(x_0^m g^h) - \deg f_j \le m+\deg(g).$$ Now in the
primary decomposition of $I^H$, suppose that for some $k$, $$x_0^{k}
g^h \in \Bigcap{i\in S}{} q_i, \mbox{ and } x_0^{k} g^h \not\in \q_j
\mbox{ if } j \not\in S.$$ Choose $\ell \not\in S$ such that
$V(q_\ell)$ is maximal.  Since $g \in I^A$, we know $V(\q_\ell)
\subset V((x_0))$, hence $x_0 \in \p_\ell$.  Let $$I_S= \Bigcap{i \in
S}{} \q_i.$$ Then $\mult_I(\q_\ell)$ is easily seen to be the length
of a maximal chain of ideals between: $$I \cdot k[\x]_{p_\ell} \mbox{
and } I_S \cdot k[\x]_{p_\ell}.$$ But look at the ideals $J_p$, for $p
\ge 0$, defined by $$I \, k[\x]_{\p_\ell} \subset
\underbrace{(I, x_0^{k+p} g^h) \, k[\x]_{\p_\ell}}_{\textstyle J_p}
\subset I_S \, k[\x]_{\p_\ell}.$$
If $J_p=J_{p+1}$, then
\begin{eqnarray*}
x_0^{k+p} g^h &\in& (I, x_0^{k+p+1} g^h) \\ \mbox{i.e.,} \;\;
x_0^{k+p} g^h &=& a \, x_0^{k+p+1} g^h + b, \;\; b \in I.
\end{eqnarray*}
But $1- a x_0$ is a unit in $k[\x]_{\p_\ell}$, so $$J_p= x_0^{k+p}
g^h= (1- a x_0)^{-1} b \in I \, k[\x]_{\p_\ell}.$$ This means that in
any case $$x_0^{k+\mult_I(\q_\ell)} g^h \in I \cdot k[\x]_{\p_\ell}$$
hence, because $q_\ell$ is $p_\ell$-primary:
$$x_0^{k+\mult_I(\q_\ell)} g^h \in \q_\ell$$ Induction now shows that
$$x_0^{\mult_\infty(I^H)} g^h \in I^H$$

\begin{corollary}\label{cor10}
Let $I^A$, $I^H$ be as above.  If $g \in I^A$, then
$$g= \sum a_i f_i$$
where $\deg(a_i) \le \deg(g)+{m(I)+n+1 \choose n+1}$.
\end{corollary}

\begin{proofn}
Combine Propositions \ref{prop4} and \ref{prop9}.
\end{proofn}

If we further estimate $m(I)$ by \thmref{thm4} in characteristic
$0$ or by \propref{prop5}, we get somewhat weaker versions of
Hermann's \thmref{thm1}.  But if $I=V(Z)$, $Z$ a good variety, we
may expect the Corollary to give much better bounds than \thmref{thm1}.

\corref{cor10} shows that any example which demonstrates the
necessity of double exponential growth in Hermann's ideal membership
bound (\thmref{thm1}) also demonstrates the necessity of double
exponential growth in the bounds on $m(I)$ given in \thmref{thm4}
and \propref{prop5}.  Thus we can make use of the general
arguments for the existence of such examples given in
\cite{mm82}, rather than depending on the single example of
\propref{prop7}, to show that the bounds on $m(I)$ inevitably grow double
exponentially: Since in \corref{cor10}, the degrees of the $a_i$ are
bounded by a single exponential function of $m(I)$, in all examples
where the degrees of the $a_i$ grow double exponentially, $m(I)$ also
grows double exponentially.

This line of argument gives a geometric link between the ideal
membership problem and $m(I)$: In \corref{cor10}, if $I^A$
exhibits $a_i$ of high degree, then $I^H$ has primary components of
high multiplicity.  These components force $m(I)$ to be large, and
distinguish $I^H$ from good ideals considered in \thmref{thm8}
and related conjectures.

A major step in understanding the gap between the double exponential
examples and the strong linear bounds on the regularity of many smooth
varieties was taken by Brownawell \cite{bro87} and Koll{\'a}r
\cite{kol88}.  They discovered the beautiful and satisfying fact that
if we replace membership in $I$ by membership in $\sqrt{I}$, then
there are single exponential bounds on the degrees of $a_i$:

\begin{thm}[Brownawell, Koll{\'a}r]
Let $k$ be any field, let $I = (f_1,...,f_k) \subset k[x_1,...,x_n]$
and let $d = \max(\deg(f_i),i=1,\cdots,k; 3)$. If $n=1$, replace $d$
by $2d-1$.  If $g \in \sqrt{I}$, then there is an expression $$ g^s =
\sum_{i=1}^k a_i f_i $$ where $s \leq d^n$ and $\deg(a_i) \le (1 +
\deg(g))d^n$.  In particular: $$
\left( \sqrt{I} \right) ^{d^n} \subset I.
$$
\end{thm}

What this shows is that although the bad examples have to have primary
components at infinity of high degree, nonetheless these primary
ideals contain relatively small powers of $\sqrt{I^H}$.  The picture
you should have is that these embedded components at infinity are like
strands of ivy that creep a long way out from the hyperplane at
infinity, but only by clinging rather closely to the affine
components.

\vskip .25in

{\bf\noindent Technical Appendix to \secref{bounds}}

\vskip .15in
\noindent
{\bf 1.}  Proof of the equivalence of the conditions in \defref{def1}:

In \cite[pp. 99-101]{mum66}, it is proven that for any coherent sheaf
$\Fh$ on $\proj n$, $H^i(\Fh(-i))=(0)$, $i \ge 1$ implies that the
same holds for $\Fh(d)$, all $d \ge 0$, and that $H^0(\Fh(d))$ is
generated by $H^0(\Fh) \otimes H^0(\Oh(d))$.  In particular, if you
apply this to $\Fh = \Ih(m)$, the equivalence of (a) and (b) follows.
(Note the diagram: $$
\begin{array}{ccc}
I_d & \longrightarrow & H^0(\Ih(d)) \\
\bigcap & & \bigcap \\
k[\x]_d & \longrightarrow & H^0(\Oh_{\proj n}(d))
\end{array}
$$ which shows that $I_m \rightarrow H^0(\Ih(m))$ is injective for
every $d$).  To show that (b) $\Rightarrow$ (c), first note that we
may rephrase the reults in \cite{mum66} to say that if
$H^i(\Fh(-i))=(0)$, $i \ge 1$, then the degrees of the minimal
generators of the $k[\x]$-module $$\stackrel{\textstyle\oplus}{_{d \in
\ZZ}} \;\; H^0(\Fh(d))$$ are all zero or less.  So we may construct
the resolution in (c) inductively: at the $k$\thuh{th} stage, say $$
\Limits{\oplus}{\alpha = 1}{r_k + 1} \; k[\x] \cdot e_{\alpha,k-1}
\stackrel{\phi_{k-1}}\longrightarrow \cdots \longrightarrow k[\x]
\longrightarrow k[\x]/I \longrightarrow 0
$$ has been constructed, let $M_k= \mbox{ ker}(\phi_{k-\ell})$ and let
$\Fh_k$ be the corresponding sheaf of ideals.  The induction
hypothesis will say that $H^i(\Fh_k(m+k-1))=(0)$, $i \ge1$.  Therefore
$M_k$ is generated by elements of degree $ \le m+k$, i.e.,
$d_\alpha=\deg e_{\alpha,k} \le m+k$, all $\alpha$.  We get an exact
sequence $$0 \;\longrightarrow\; M_{k+1} \;\longrightarrow
\Limits{\oplus}{\alpha=1}{r_k} \left( k[\x] \cdot e_{\alpha,k} \right)
\;\longrightarrow\; M_k \;\longrightarrow\; 0
$$
hence
\begin{eqnarray}
0 \;\longrightarrow\; \Fh_{k+1} \;\longrightarrow
\Limits{\oplus}{\alpha=1}{r_k} \Oh_{\proj n}(-d_\alpha)
\;\longrightarrow\; \Fh_k \;\longrightarrow\; 0
\end{eqnarray}
Therefore
\begin{eqnarray}
\Limits{\oplus}{\alpha=1}{r_k} H^i(\Oh_{\proj
n}(m\!+\!k\!-\!i\!-\!d_\alpha))
\longrightarrow H^i(\Fh_k(m\!+\!k\!-\!i)) \longrightarrow \hskip 0.3in \\
\;\;\;\;H^{i+1}(\Fh_{k+1}(m\!+\!(k\!+\!1)\!-\!(i\!+\!1)))
\longrightarrow \Limits{\oplus}{\alpha=1}{r_k} H^{i+1}(\Oh_{\proj
n}(m\!+\!k\!-\!i\!-\!d_\alpha)) \nonumber
\end{eqnarray}
is exact.  But $m+k-i-d_\alpha \ge -i$ so $H^{i+1}(\Oh_{\proj
n}(m+k-i-d_\alpha))=(0)$.  This shows that $\Fh_{k+1}$ satisfies the
induction hypothesis and we can continue.  Thus (c) holds.  To see
that (c) $\Rightarrow$ (a), we just use the same exact sequences (1)
and prove now by descending induction on $k$ that
$H^i(\Fh_k(m+k-i))=(0)$, $i \ge1$.  Since $I=\Fh_0$, this does it.
The inductive step again uses (2), since $H^i(\Oh_{\proj
n}(m+k-i-d_\alpha))=(0)$ too.

\vskip .15in
\noindent
{\bf 2.}  Proof of \propref{prop4}:

Look first at the case $r=0$.  Let $\Ih$ be the sheaf of ideals
defined by $I$ and let $\Ih^* \supset \Ih$ be the sheaf defined by
omitting all $0$-dimensional primary components of $I$.  Consider the
exact sequence: $$0 \;\longrightarrow\; \Ih(m-1) \;\longrightarrow
\Ih^*(m-1) \;\longrightarrow\; (\Ih^*/\Ih)(m-1) \;\longrightarrow\; 0
$$
This gives us:
$$
H^0(\Ih^*(m-1)) \;\longrightarrow\; H^0((\Ih^*/\Ih)(m-1))
\;\longrightarrow\; H^1(\Ih(m-1))
$$ Now $H^1(\Ih(m-1))=(0)$ by $m$-regularity, and
$h^0((\Ih^*/\Ih)(m-1)) = h^0(\Ih^*/\Ih) = \mbox{ length}(\Ih^*/\Ih) =
\arithdeg_0(I)$ since $\Ih^*/\Ih$ has $0$-dimensional support. But
$H^0(\Ih^*(m-1)) \subset H^0(\Oh_{\proj n}(m-1))$, so
\begin{eqnarray*}
\arithdeg_0(I) & \le & h^0(\Ih^*(m-1)) \\
& \le & h^0(\Oh_{\proj n}(m-1)) \\
& = & { m+ n -1 \choose n }
\end{eqnarray*}

If $r>0$, we can prove the Proposition by induction on $r$.  Let $H$
be a generic hyperplance in $\proj n$, given by $h=0$.  Let $I_H=
(I,h) / (h) \subset k[x_0,\ldots,x_n]/(h) \cong
k[x_0^\prime,\ldots,x_{n-1}^\prime]$ for suitable linear combinations
$x_i^\prime$ of $x_i$.  Then it is easy to check that:
$$\arithdeg_r(I) = \arithdeg_{r-1}(I_H)$$ and that $I_H$ is also
$m$-regular, so by induction
\begin{eqnarray*}
\arithdeg_{r-1}(I_H) & \le & {m+(n-1)-(r-1)-1 \choose (n-1)-(r-1)} \\
& = & { m+n-r-1 \choose n-r }
\end{eqnarray*}
If $r = -1$, we use the fact that $$0 \longrightarrow I_d
\longrightarrow H^0(\Ih(d)) \stackrel{\approx}{\longleftarrow}
(I^{\mbox{sat}})_d$$ if $d \ge m$, hence $$\dim(I^{\mbox{sat}}/I) \le
\dim k[\x]/(x_0,\ldots,x_n)^m = {m+n
\choose n+1}.  $$

\vskip .15in
\noindent
{\bf 3.}  Proof of \propref{prop5}:

Let $I \subset k[x_0, \ldots x_n]$ and assume, after a linear change
of coordinates, that $x_n$ is not contained in any associated prime
ideals of $I$.  Let $\overline I \subset k[x_0, \ldots x_{n-1}]$ be
the image of $I$.  Then $d(\overline I) = d(I)$ and by induction we
may assume $$m(\overline I) \le (2d(I))^{(n-1)!}.$$ We will prove, in
fact, that
\begin{eqnarray}
m(I) \le m(\overline I) + { m(\overline I) - 1 + n \choose n }
\end{eqnarray}
and then we will be done by virtue of the elementary estimate:
$$\mbox{if } m^* = (2d(I))^{(n-1)!}, \mbox{ and } d \ge 2,
\mbox{ then } m^* + { m^* - 1 + n \choose n } \le  (2d(I))^{n!}$$
To prove (3), we use the long exact sequence
$$
\begin{array}{ccccccccc}
0 & \longrightarrow & (I:(x_0))_{k-1} & \stackrel{x_0}\longrightarrow
& I_k & \longrightarrow & \overline I_k & \longrightarrow & 0 \\ & &
\downarrow & & \downarrow & & \downarrow \\ 0 & \longrightarrow &
H^0(\Ih(k-1)) & \longrightarrow & H^0(\Ih(k-1)) & \longrightarrow &
H^0(\overline\Ih(k-1)) & \stackrel{\delta}\longrightarrow \\ \\ &
\stackrel{\delta}\longrightarrow & H^1(\Ih(k-1)) & \longrightarrow &
H^1(\Ih(k)) & \longrightarrow & H^1(\overline\Ih(k)) & \\
\end{array}
$$ where $(I:(x_0)) = \setdef{f}{x_0 f \in I}$.  Let $\overline m =
m(\overline I)$.  Note that $H^i(\overline\Ih(k-1)) = (0)$, $i \ge 1$,
$k \ge \overline m$, hence $$H^i(\Ih(k-1)) \rightarrow H^i(\Ih(k))$$
is an isomorphism if $k \ge \overline m-1+1$ and $i \ge 2$.  Since
$H^i(\Ih(k)) = (0)$, $k \gg 0$, this shows that $H^i(\Ih(k)) = (0)$,
$i \ge 2$, $k \ge \overline m -i$.  Moreover $\overline I_k
\rightarrow H^0(\overline\Ih(k))$ is an isomorphism if $k \ge
\overline m$, hence $\delta = 0$ if $k \ge \overline m$, hence
$H^1(\Ih(k)) = (0)$, $k \ge
\overline m - 1$.  But now look at the surjectivity of $I_k \rightarrow
H^0(\Ih(k))$.  For all $k$, let $M_k$ be the cokernel.  Then $\oplus
M_k$ is a $k[\x]$-module of finite dimension.  Multiplication by $x_0$
induces a sequence: $$0 \; \longrightarrow \; { (I:(x_0))_{k-1} \over
I_{k-1} }
\; \longrightarrow \; M_{k-1}
\; \stackrel{x_0}\longrightarrow \; M_k \; \longrightarrow \; 0$$
which is exact if $k \ge \overline m$.  But if, for one value of $k
\ge \overline m$, \begin{eqnarray} (I:(x_0))_k = I_k
\end{eqnarray}
then by \thmref{thm2}, $I$ is $k$-regular and (4) continues to hold
for larger $k$, and $M_k$ must be $(0)$.  In other words, $$\dim M_k,
\;\; k \;\ge\; \overline m-1$$ is non increasing and monotone
decreasing to zero when $k \ge
\overline m$.  Therefore
\begin{eqnarray*}
m(I) & \le & \overline m + \dim M_{\overline m-1} \\
& \le & \overline m + \dim  k[\x]_{\overline m-1} \\
& \le & \overline m + {\overline m-1+n \choose n} \\
\end{eqnarray*}
which proves (3).

\vskip .15in
\noindent
{\bf 4.}   Proof of \thmref{thm8}(b):

Let $Z$ be a smooth $r$-dimensional subvariety of $\proj n$ and $d =$
degree of $Z$.  We first consider linear projections of $Z$ to $\proj
r$ and to $\proj{r-1}$.  To get there, let $L_1 \subset \proj n$ be a
linear subspace of dimension $n-r-1$ disjoint from $Z$ and $L_2
\subset L_1$ a linear subspace of dimension $n-r-2$.  Take these as
centers of projection: $$
\begin{array}{ccc}
\proj{n} - L_1 & \supset & \;\;\;\; \;Z \;\;\;
\stackrel{p_2}\longrightarrow
\\ \downarrow & & \downarrow p_1 \\
\proj{r+1} - \{ P \} & \supset & Z_1 \;\;\; \\
\downarrow \\
\stackrel{p_2}\longrightarrow \; \proj{r} \;\;\;\;\;
\end{array}
$$ Let $x_0, \ldots x_{r+1}$ be coordinates on $\proj{r+1}$ so that $p
= (0, \ldots , 0, 1)$, hence $x_0, \ldots x_r$ are coordinates on
$\proj r$.  Let $f(x_0, \ldots x_{r+1}) = 0$ be the equation of the
hypersurface $Z_1$.

Now there are two ways of getting $r$-forms on $Z$: by pullback of
$r$-forms on $\proj r$ and by residues of $(r+1)$-forms on
$\proj{r-1}$ with simple poles along $Z_1$.  The first gives us a
sheaf map $$p_2^* \;\Omega_{\proj r}^r \hookrightarrow \Omega_Z^r$$
whose image is $\Omega_Z^r(-B_1)$, $B_1$ the branch locus of $p_2$.
Corresponding to this on divisor classes:
\begin{eqnarray}
K_Z & \equiv & p_2^*(K_{\proj r}) + B_1 \\
& \equiv & - (r+1)H + B_1, \nonumber
\end{eqnarray}
where $H =$ hyperplane divisor class on $Z$.  The
second is defined by
\begin{eqnarray}
a(\x) \cdot { dx_1 \wedge \ldots \wedge dx_{r+1} \over f }
\;\longmapsto\; p_1^* \left(a(\x) \cdot { dx_1 \wedge \ldots \wedge
dx_r \over
\partial f/\partial x_{r+1}} \right)
\end{eqnarray}
and it gives us an isomorphism $$p_1^*(\Omega_{\proj{r+1}}^{r+1}(Z_1)
{\Large\mid}_{Z_1} )
\;\cong\; \Omega_Z^r(B_2)$$
$B_2$ is a divisor which can be interpreted as the {\em conductor} of
the affine rings of $Z$ over those of $Z_1$: i.e., $$f \in \Oh_Z(-B_2)
\;\Longleftrightarrow\; f \cdot (p_{1,*} \Oh_Z) \subset
\Oh_{Z_1}.$$
In particular,
\begin{eqnarray}
p_{1,*}(\Oh_Z(-B_2)) \cong \mbox{ sheaf of } \Oh_{Z_1} - \mbox{ ideals
} C
\mbox{ in } \Oh_{Z_1}.
\end{eqnarray}

A classical reference for these basic facts is Zariski \cite{zar69},
Prop.  12.13 and Theorem 15.3.  A modern reference is Lipman
\cite{lip84} (apply Def. (2.1)b to $p_1$ and apply Cor. (13.6) to $Z_1
\subset \proj{r+1}$).  (4) gives us the divisor class identity:
\begin{eqnarray}
K_Z + B_2 & \equiv & p_1^*(K_{\proj{r+1}} + Z_1) \\
& \equiv & (d-r-2) H. \nonumber
\end{eqnarray}
(5) and (8) together tell us that $$B_1 + B_2 \equiv (d-1)H.$$ In
fact, the explicit description (6) of the residue tells us more:
namely that if $y_1,\ldots,y_r$ are local coordinates on $Z$, then
$${\partial (x_1,\ldots,x_r) \over \partial (y_1,\ldots,y_r)} \cdot {1
\over \partial f / \partial x_{r+1}} \; dy_1 \wedge \ldots \wedge
dy_r$$ generates $\Omega_Z^r(B_2)$ locally.  But ${\partial
(x_1,\ldots,x_r) \over \partial (y_1,\ldots,y_r)} = 0$ is a local
equation for $B_1$, so this means that $\partial f / \partial x_{r+1}
= 0$ is a local equation for $B_1+B_2$.  But $\partial f / \partial
x_{r+1} = 0$ is a global hypersurface of degree $d-1$ in $\proj{r+1}$,
hence globally: $$B_1 + B_2 = p_1^*(V({\partial f \over \partial
x_{r_1}}))$$ (equality of divisors, not merely divisor classes).  All
this is standard classical material.

(7) has an important cohomological consequence: let $C^* \subset
\Oh_{\proj{r+1}}$ be the sheaf of ideals consisting of functions whose
restriction to $Z_1$ lies in C.  Then we get an exact sequence: $$0
\,\rightarrow\,
\Oh_{\proj{r+1}}(-Z_1) \,\rightarrow\,
C^* \, \Oh_{\proj{r+1}} \,\rightarrow\,
C \, \Oh_{Z_1} \,\rightarrow\, 0$$
hence an exact sequence
$$0 \,\rightarrow\,
\Oh_{\proj{r+1}}(\ell-d) \,\rightarrow\,
C^* \Oh_{\proj{r+1}}(\ell) \,\rightarrow\, p_{1,*}(\Oh_Z(\ell H -
B_2)) \,\rightarrow\, 0$$ for all integers $\ell$.  But
$H^1(\Oh_{\proj{r+1}}(\ell-d)) = (0)$, hence $$H^0(C^*
\Oh_{\proj{r+1}}(\ell)) \,\rightarrow\, H^0(\Oh_Z(\ell H - B_2))$$ is
surjective, hence
\begin{eqnarray}
H^0(\Oh_Z(\ell H - B_2)) \,\subset\, \mbox{ Im}\left[
H^0(\Oh_{\proj{r+1}}(\ell))
\,\rightarrow\, H^0(\Oh_Z(\ell H)) \right].
\end{eqnarray}

Now let us vary the projections $p_1$ and $p_2$.  For each choice of
$L_1$, we get a different $B_1$: call it $B_1(L_1)$, and for each
choice of $L_2$, as different $B_2$: call it $B_2(L_2)$.  By (5) and
(8), all divisors $B_1(L_1)$ are linearly equivalent as are all
divisors $B_2(L_2)$.  Moreover:
\begin{eqnarray*}
\Bigcap{L_1}{} B_1(L_1) &=& \emptyset \\
\Bigcap{L_2}{} B_2(L_2) &=& \emptyset
\end{eqnarray*}
This is because, if $x \in Z$, then there is a choice of $L_1$ such
that $p_1: Z \rightarrow \proj{r}$ is unramified at $y$; and a choice
of $L_2$ such that $p_2(x) \in Z_1$ is smooth, hence $p_2$ is an
isomorphism near $x$.  Thus $$\abs{B_1(L_1)} = \abs{K_Z + (r+1) H}$$
and $$\abs{B_2(L_2)} = \abs{K_Z + (d-r-2) H}$$ are base point free
linear systems.

Next choose $(r+1)$ $L_2$'s, called $L_2^\alpha$, $1 \le \alpha \le
r+1$, so that if $B_2^{(\alpha)} = B_2(L_2^{(\alpha)})$, then
$\Bigcap{\alpha}{} B_2^{(\alpha)} =
\emptyset$.  Look at the Koszul complex:
\begin{eqnarray*}
0 \,\rightarrow\, \Oh_Z(\ell H - \sum B_2^{(\alpha)}) &\rightarrow&
\cdots \\
\rightarrow\, \sum_{\alpha,\beta} \Oh_Z(\ell H -
B_2^{(\alpha)}-B_2^{(\beta)}) &\rightarrow& \sum_{\alpha}\Oh_Z(\ell H
- B_2^{(\alpha)}) \,\rightarrow\, \Oh_Z(\ell H) \,\rightarrow\, 0.
\end{eqnarray*}

This is exact and diagram chasing gives the conclusion: $$
H^i(\Oh_Z(\ell H - (i+1) B_2)) = (0), \mbox{ all } i \ge 1 $$
$$\Rightarrow \sum_\alpha H^0(\Oh_Z(\ell H - B_2^{(\alpha)}))
\rightarrow H^0(\Oh_Z(\ell H)) \mbox{ surjective} $$ hence by (9)
$$H^0(\Oh_{\proj{n}}(\ell)) \rightarrow H^0(\Oh_Z(\ell H)) \mbox{
surjective}$$ and $$H^{i+j}(\Oh_Z(\ell H - (i+1)B_2)) = (0), \mbox{
all } i \ge 0$$ $$\Rightarrow \;\; H^j(\Oh_Z(\ell H) = (0).$$ Now
$I(Z)$ is $m$-regular if and only if $H^i(\Ih_Z(m-i)) = (0)$, $i \ge
1$, hence if and only if $$H^0(\Oh_{\proj{n}}(m-1)) \rightarrow
H^0(\Oh_Z(m-1)) \mbox{ surjective}$$ $$H^i(\Oh_Z(m-i-1)) = (0), \; i
\ge 1. $$ By the previous remark, this follows provided that
$$H^{i+j}(\Oh_Z((m-i-1)H - (j+1)B_2)) = (0), \mbox{ if } i,j \ge 0,\;
i+j \ge 1.$$ But let us rewrite: $$(m-i-1)H-(j+1)B_2 \;\equiv\; K_Z
+jB_1+(m-i-(j+1)(d-1)+r)H$$ using (5) and (8).  Note that $j B_1 +
\ell H$ is an ample divisor if $\ell \ge 1, j \ge 0$, because
$\abs{B_1}$ is base point free.  Therefore by the Kodaira Vanishing
Theorem, $$H^i(\Oh_Z (K_Z+j B_1 + \ell H)) = (0), \; i,j \ge 1, \; j
\ge 0$$ and provided $m=(r+1)(d-2)+2$, this gives the required
vanishing.

\section{Applications}
\label{applications}
{}From some points of view, the first main problem of algebraic geometry
is to reduce the study of a general ideal $I$ to that of prime ideals,
or the study of arbitrary schemes to that of varieties. One way of
doing this is to find a decomposition of the ideal into primary
ideals: i.e. write it as an intersection of primary ideals. But even
when non-redundancy is added, this is not unique, and usually one
actually wants something less: to find its radical and perhaps write
the radical as an intersection of prime ideals, or to find its top
dimensional part, or to find its associated prime ideals and their
multiplicities.  There are really four computational problems involved
here which should be treated separately: (i) eliminating the
multiplicities in the ideal $I$, (ii) separating the pieces of
different dimension, (iii) ``factoring'' the pieces of each dimension
into irreducible components, and finally (iv) describing the original
multiplicities, either numerically or by a primary ideal. Three of
these four problems are the direct generalizations of the basic
problems for factoring a single polynomial: we can eliminate multiple
factors, getting a square-free polynomial, we can factor this into
irreducible pieces and we can ask for the multiplicities with which
each factor appeared in the original polynomial. There is a fifth
question which arises when we work, as we always must do on a
computer, over a non-algebarically closed field $k$: we can ask (v)
for an extension field $k'$ of $k$ over which the irreducible
components break up into absolutely irreducible components.

Classical algorithms for all of these of these rely heavily on making explicit
projections of $V(I)$ to lower dimensional projective spaces. This can be done
either by multi-variable resultants if you want only the set-theoretic
projection, or by \std\ bases with respect to the lexicographic order or an
elimination order, to get the full ideal $I \cap k[X_0, \cdots, X_m]$. Recent
treatments of multi-variable resultants can be found in \cite{can89},
\cite{cha91}, and a recent treatment of the basis method can be found in
\cite{gtz88}. There is no evidence that either of these is an efficient method,
however, and taking \std\ bases in the lexicographical order or an elimination
order is often quite slow, certainly slow in the worst case. The general
experience is that taking projections can be very time consuming. One reason is
that the degree of the generators may go up substantially and that sparse
defining polynomials may be replaced by more or less generic polynomials. A
specific example is given by principally polarized abelian varieties of
dimension $r$: they are defined by quadratic polynomials in $(4^r - 1)$-space,
but their degree here (hence the degree of their generic projection to
$\proj{r+1}$) is $4^r r!$ \cite{mum70a}. In fact, any variety is defined
purely by quadratic relations in a suitable embedding \cite{mum70b}.

Instead of using real computational experience, the fundamental method in
theoretical computer science for analyzing complexity of algorithms is to count
operations. For algebraic algorithms, the natural measure of complexity is not
the number of bit operations, but the number of field operations, addition,
subtraction, multiplication and (possibly) division that are used. In this
sense, any methods that involve taking \std\ bases for any order on monomials
will have a worst-case behavior whose complexity goes up with the regularity of
the ideal hence will take ``double exponential time''. However, it appears that
this worst-case behavior may in fact only concern problem (iv) -- finding the
primary ideals -- and that problems (i), (ii) and (iii) may be solvable in
``single exponential time''. The idea that such algorithms should exist for
finding $V(I)$ set-theoretically was proposed in the 1984 lecture on which this
article is based, but turned out, in fact, to have been already proven by
Chistov and Grigoriev, cf.  their unpublished 1983 note \cite{cg83}. Their line
of research led, in some sense, to the work of Brownawell and Koll{\'a}r,
showing the single exponential bound $\left(\sqrt{I}\right)^m \subset I$ for $m
= d^n$, where $d = \max(\hbox{ degrees of generators of }I)$.

Based on this work, Giusti and Heintz \cite{gh91} give a singly
exponentially bounded algorithm for computing ideals ${\bf q}_i$ such
that $V({\bf q}_i)$ are the irreducible components of $V(I)$ (over the
ground field $k$). The method depends on computing what is essentially
the Chow form of each component, and leads to an ideal defining this
variety but not its full ideal. In fact, their ${\bf q}_i$ may be
guaranteed to be prime except for possible embedded components.

A direct approach to constructing both $\sqrt{I}$ and the intersection
of the top-dimensional primary components of $I$, denoted
$\hbox{Top}(I)$, is given in a recent paper by Eisenbud, Huneke and
Vasconcelos \cite{ehv92}. Their construction of the radical uses the
Jacobian ideals, i.e. the ideals of minors of various sizes of the
Jacobian matrix of generators of $I$. This is certainly the most
direct approach, but, again they have trouble with possible embedded
components, and must resort to ideal quotients, hence they need a
\std\ basis of $I$ in the reverse lexicographic order. They compute
$\hbox{Top}(I)$ as the annihilator of $\hbox{Ext}^{\codim(I)}
(k[X_0,\cdots,X_n]/I, k[X_0,\cdots,X_n])$, which is readily
found from a full resolution using \std\ bases. Their algorithm
appears to be practical in some cases of interest, but still has
double exponential time worst-case behavior.

It may turn out to be most effective in practice to combine these ideas.
Often an ideal under study has regularity far smaller than the
geometric degree of its top dimensional components; projecting these
components to a hypersurface requires computing in degrees up to the
geometric degree, which is wasteful. On the other hand, methods such
as those in \cite{ehv92} work better in low codimensions, if only
because there are fewer minors to consider in the Jacobian matrix.
Thus, projecting an arbitary scheme down to low codimension and then
switching to direct methods may work best of all.

This still does not settle the issue of the complexity of calculating
$\sqrt{I}$, or, for that matter, calculating the full prime ideal of
any subvariety of codimension greater than one. Chow form type methods
give you an effective method of defining the set $V(I)$ but only of
generating $I$ up to possible embedded components. For this reason,
the two schools of research, one based on the algebra of $I$, the
other based on subsets of ${\bf P}^n$ have diverged. If we knew, as
discussed in the previous section, that the regularity of a reduced
ideal could be bounded singly exponentially, then we could bound the
degrees of the generators of $\sqrt{I}$, and, using Brownawell-Koll{\'a}r,
we could determine $\sqrt{I}$ up to these degrees and get the whole
ideal. But without such a bound, it is still not clear whether only
$V(I)$ and not $\sqrt{I}$ can be found in worst-case single
exponential time.

Let's look at problem (iii). Assume you have found a reduced
equidimensional $I$. To study splitting it into irreducible or
absolutely irreducible pieces, we shall assume initially it is a
hypersurface, i.e.  $I=(f)$. Computationally, there may often be
advantages to not projecting a general $I$ to a hypersurface, and we
will discuss one such approach below. Geometrically, there is nothing
very natural about irreducible but not absolutely irreducible
varieties: from the standpoint of their properties, they behave like
reducible varieties, except that, being conjugate over $k$, their
components have very similar properties. If the ground field $k$ gets
bigger or smaller, the set of absolutely irreducible components gets
partitioned in finer or coarser ways into the $k$-components. If one
has never done any calculations, one would therefore be inclined to
say -- let's extend $k$ as far as needed to split our algebraic set up
into absolutely irreducible components. {\it This is a very bad idea!}
Unless this extension $k'$ happens to be something simple like a
quadratic or cyclotomic extension of $k$, the splitting field $k'$ is
usually gigantic. This is what happens if one component of $V(I)$ is
defined over an extension field $k_1$ of $k$ of degree $e$, and the
Galois group of $k_1/k$ is the full symmetric group, a very common
occurence. Then $V(I)$ only splits completely over the Galois closure
of $k_1/k$ and this has degree $e!$. The moral is: never factor unless
you have to.

In fact, unless you need to deal simultaneously with more than one of
its irreducible components, you can proceed as follows: the function
field $K = k[X_0,\cdots,X_n]/(f)$ contains as a subfield an isomorphic
copy of $k_1$: you find that field as an extension $k_1 = k[y]/(p(y))$
of $k$, and solve for the equation of one irreducible component $f_1
\in k_1[X_0,\cdots,X_n]$ by the formula $\hbox{Norm}_{k_1/k}(f_1) = f$.

Pursuing this point, why should one even factor the defining equation
$f$ over $k$? Factoring, although it takes polynomial time
\cite{lll82}, is often very slow in real time, and, unless the
geometry dictates that the components be treated separately, why not
leave them alone. In some situations, for instance, \cite{dd84} one
may have an ideal, module or other algebraic structure defined by
polynomials or matrices of polynomials over a {\it ground ring} $D =
k[y]/(p(y))$, where $p$ is a square-free polynomial. Thus $D$ is a
direct sum of extension fields, but there is no need to factor $p$ or
split up $D$ until the calculations take different turns with the
structures over different pieces of $Spec(D)$.

The standard methods of factoring in computer algebra all depend on
(i) writing the polynomial over a ring, finitely generated over $\bf
Z$, and reducing modulo a maximal ideal $\bf m$ in that ring,
obtaining a polynomial over a finite field; and (ii) restricting to a
line $L$, i.e.  substituting $X_i = a_i X_0 + b_i, i \geq 1$ for all
but one variable, obtaining a polynomial in one variable over a finite
field. This is then factored and then, using Hensel's lemma, one lifts
this factorization modulo higher powers of $\bf m$ and of the linear
space $L$. One then checks whether a coarsened version of this
factorization works for $f$.  This is all really the arithmetic of
various small fields. Geometrically, every polynomial in one variable
factors over a suitable extension field and the question of counting
the absolutely irreducible components of a variety is really more
elementary: it is fundamentally topological and not arithmetic. One
should, therefore, expect there to be direct geometric ways of
counting these components and separating them.  Assuming $I$ is a
reduced, equi-r-dimensional ideal, the direct way should be to use
Serre duality, computing the cohomology $H^r(\Omega^r_{V(I)})$, where
$\Omega^r_{V(I)} \subset \omega_{V(I)}$ is the subsheaf of the
top-dimensional dualizing sheaf of $V(I)$ of absolutely regular
$r$-forms. Its dimension will be the number of absolutely irreducible
components into which $V(I)$ splits. Calculating this cohomology
involves two things: algebraically resolving the ideal $I$ and
geometrically resolving the singularities of $V(I)$ far enough to work
out $\Omega^r_{V(I)}$. Classically, when $I = (f)$ was principal,
$\Omega^r_{V(I)}$ was called its ideal of ``subadjoint'' polynomials.

There is one case where this is quite elementary and has been carried
out: this is for plane curves. One can see immediately what is
happening by remarking that a non-singular plane curve is
automatically absolutely irreducible, hence one should expect that its
singularities control its decomposition into absolutely irreducible
pieces. Indeed, if ${\bf C}
\subset k[X_0,X_1,X_2]/(f)$ is the conductor ideal, then
$\Omega^1_{V(f)}$ is given by the homogeneous ideal ${\bf C}$, but
with degree 0 being shifted to be polynomials of degree $d-3$, $d$ the
degree of $f$. To calculate $H^1$, assume $X_0$ is not zero at any
singularity of $V(f)$ and look at the finite-dimensional vector space
of all functions $k[X_1/X_0,X_2/X_0]/({\bf C} + (f))$ modulo the
restrictions $g/(X_0^{d-3})$ for all homogeneous polynomials $g$ of
degree $d-3$.  This will be canonically the space of functions on the
set of components of $V(f)$ with sum $0$. In particular, it is $(0)$
if and only if $V(f)$ is absolutely irreducible. This follows from
standard exact sequences and duality theory. It was known classically
as the Cayley-Bacharach theorem, for the special case where $V(f)$ was
smooth except for a finite number of ordinary double points. It states
that $V(f)$ is absolutely irreducible if and only if for every double
point $P$, there is a curve of degree $d-3$ passing through all the
double points except $P$.

This example gives one instance where a deeper computational analysis of
varieties requires a computation of its resolution of singularities. We believe
that there will be many instances where practical problems will require such an
analysis. In many ways, resolution theorems look quite algorithmic, and, for
instance, Abhyankar and his school have been approaching the problem in this
way \cite{abh82}, as have Bierstone and Milman \cite{bm91}. However, the only
case of resolution of singularities to be fully analyzed in the sense of
computational complexity is that of plane curves. This has been done by
Teitelbaum \cite{tei89}, \cite{tei90}. His analysis is notable in various ways:
he is extremely careful about not making unnecessary factorizations, let alone
taking unnecessary field extensions, and uses the ``$D$'' formalism discussed
above. He describes his algorithm so precisely that it would be trivial to
convert it to code and, as a result, he gives excellent bounds on its
complexity.


\begin{thebibliography}{Mum70b}

\bibitem[Abh82]{abh82}
S.~S. Abhyankar, {\em Weighted expansions for canonical desingularization},
  Lecture Notes in Math., vol. 910, Springer-Verlag, 1982.

\bibitem[Art76]{art76}
M.~Artin, {\em Lectures on deformations of singularities}, Tata Institute on
  Fundamental Research, Bombay, 1976.

\bibitem[Bay82]{bay82}
Dave Bayer, {\em The division algorithm and the {H}ilbert scheme}, {Ph.D.}
  thesis, Harvard University, Department of Mathematics, June 1982, order
  number 82-22588, University Microfilms International, 300 N. Zeeb Rd., Ann
  Arbor, MI 48106.

\bibitem[BEL91]{bel91}
Aaron Bertram, Lawrence Ein, and Robert Lazarsfeld, {\em Vanishing theorems, a
  theorem of {S}everi, and the equations defining projective varieties}, J.
  Amer. Math. Soc. {\bf 4} (1991), 587--602.

\bibitem[Ber78]{ber78}
G.~M. Bergman, {\em The diamond lemma for ring theory}, Adv. in Math. {\bf 29}
  (1978), 178--218.

\bibitem[BM88]{bm88}
Dave Bayer and Ian Morrison, {\em Standard bases and geometric invariant theory
  {I}. {I}nitial ideals and state polytopes}, J. Symb. Comput. {\bf 6} (1988),
  no.~2--3, 209--217, reprinted in \cite{rob89}.

\bibitem[BM91]{bm91}
E.~Bierstone and P.~Milman, {\em A simple constructive proof of canonical
  resolution of singularities}, Effective methods in algebraic geometry
  (Castiglioncello, 1990), Progr. Math., vol.~94, Birkhauser Boston, 1991,
  pp.~11--30.

\bibitem[BCR91]{bcr91}
A.~M. Bigatti, M.~Caboara, and L.~Robbiano, {\em On the computation of
  {H}ilbert-{P}oincare series}, Applicable Algebra in Engineering,
  Communications, and Computing {\bf 2} (1991), 21--33.

\bibitem[Bri73]{bri73}
J.~Briancon, {\em {W}eierstrass prepare a la {H}ironaka}, Ast{\'e}risque {\bf
  7,8} (1973), 67--73.

\bibitem[Bro87]{bro87}
W.~D. Brownawell, {\em Bounds for the degrees in the {N}ullstellensatz}, Ann.
  of Math. (2) {\bf 126} (1987), 577--591.

\bibitem[BS87a]{bs87a}
Dave Bayer and Mike Stillman, {\em A criterion for detecting {$m$}-regularity},
  Invent. Math. {\bf 87} (1987), 1--11.

\bibitem[BS87b]{bs87b}
Dave Bayer and Mike Stillman, {\em A theorem on refining division orders by the
  reverse lexicographic order}, Duke Math. J. {\bf 55} (1987), no.~2, 321--328.

\bibitem[BS88]{bs88}
Dave Bayer and Mike Stillman, {\em On the complexity of computing syzygies}, J.
  Symb. Comput. {\bf 6} (1988), 135--147.

\bibitem[BS92a]{macaulay}
Dave Bayer and Mike Stillman, {\em Macaulay: A system for computation in
  algebraic geometry and commutative algebra}, 1982--1992, computer software
  available via anonymous ftp from {\tt \small zariski.harvard.edu}.

\bibitem[BS92b]{bs92}
Dave Bayer and Mike Stillman, {\em Computation of {H}ilbert functions}, J.
  Symb. Comput. {\bf 6} (1992), 31-50.

\bibitem[Buc65]{buc65}
B.~Buchberger, {Ph.D.} thesis, Univ. Innsbr{\"u}ck, 1965.

\bibitem[Buc76]{buc76}
B.~Buchberger, {\em A theoretical basis for the reduction of polynomials to
  canonical forms}, ACM SIGSAM Bull. {\bf 39} (1976), 19--29.

\bibitem[Buc79]{buc79}
B.~Buchberger, {\em A criterion for detecting unnecessary reductions in the
  construction of {Gr\"obner} bases}, Symbolic and Algebraic Computation
  (Proceedings of {EUROSAM} 79), Lecture Notes in Computer Science, vol.~72,
  Springer-Verlag, 1979, pp.~3--21.

\bibitem[Can89]{can89}
J.~Canny, {\em Generalized characteristic polynomials}, Symbolic and Algebraic
  Computation (Proceedings of {ISSAC} 88), Lecture Notes in Computer Science,
  vol. 358, Springer-Verlag, 1989, pp.~293--299.

\bibitem[CG83]{cg83}
A.~L. Chistov and D.~Yu. Grigoriev, {\em Subexponential-time solving systems of
  algebraic equations {I}, {II}}, {S}teklov Mathematical Institute, Leningrad
  department, LOMI Preprints E-9-93, 0E-10-c83, 1983.

\bibitem[Cha91]{cha91}
Marc Chardin, {\em Un algorithme pour le calcu des r{\'e}sultants}, Effective
  methods in algebraic geometry (Castiglioncello, 1990), Progr. Math., vol.~94,
  Birkhauser Boston, 1991, pp.~47--62.

\bibitem[DD84]{dd84}
C.~Dicrescenzo and D.~Duval, {\em Computations on curves}, Lecture Notes in
  Computer Science, vol. 174, Springer-Verlag, 1984.

\bibitem[EG84]{eg84}
David Eisenbud and Shiro Goto, {\em Linear free resolutions and minimal
  multiplicity}, J. Algebra {\bf 88} (1984), no.~1, 89--133.

\bibitem[EHV92]{ehv92}
David Eisenbud, Craig Huneke, and Wolmer Vasconcelos, {\em Direct methods for
  primary decomposition}, Invent. Math. (1992), to appear.

\bibitem[Eis92]{eis92}
David Eisenbud, {\em Commutative algebra with a view toward algebraic
  geometry}, 1992, in preparation.

\bibitem[Gal74]{gal74}
A.~Galligo, {\em A propos du theoreme de preparation de {W}eierstrass},
  Fonctions de Plusieurs Variables Complexes, Lecture Notes in Math., vol. 409,
  Springer-Verlag, 1974, pp.~543--579.

\bibitem[Gal79]{gal79}
A.~Galligo, {\em Theoreme de division et stabilite en geometrie analytique
  locale}, Ann. Inst. Fourier (Grenoble) {\bf 29} (1979), 107--184.

\bibitem[GH91]{gh91}
Marc Giusti and Joos Heintz, {\em Algorithmes---disons rapides---pour la
  decomposition d'une variete algebrique en composantes irreductibles et
  equidimensionnelles [``{F}ast'' algorithms for the decomposition of an
  algebraic variety into irreducible and equidimensional components]},
  Effective methods in algebraic geometry (Castiglioncello, 1990), Progr.
  Math., vol.~94, Birkhauser Boston, 1991, pp.~169--194.

\bibitem[Giu84]{giu84}
Marc Giusti, {\em Some effectivity problems in polynomial ideal theory},
  {EUROSAM} 84), Lecture Notes in Computer Science, vol. 204, Springer-Verlag,
  1984, pp.~159--171.

\bibitem[GLP83]{glp83}
L.~Gruson, R.~Lazarsfeld, and C.~Peskine, {\em On a theorem of {C}astelnuovo,
  and the equations defining space curves}, Invent. Math. {\bf 72} (1983),
  491--506.

\bibitem[GTZ88]{gtz88}
P.~Gianni, B.~Trager, and G.~Zacharias, {\em \std\ bases and primary
  decomposition of polynomial ideals}, J. Symb. Comput. {\bf 6} (1988),
  no.~2--3, 149--167, reprinted in \cite{rob89}.

\bibitem[Her26]{her26}
Grete Hermann, {\em Die {F}rage der endlich vielen {S}chritte in der {T}heorie
  der {P}olynomideale}, Math. Ann. {\bf 95} (1926), 736--788.

\bibitem[Hir64]{hir64}
H.~Hironaka, {\em Resolution of singularities of an algebraic variety over a
  field of characteristic zero: {I}, {II}}, Ann. of Math. (2) {\bf 79} (1964),
  109--326.

\bibitem[Kol88]{kol88}
J{\'a}nos Koll{\'a}r, {\em Sharp effective {N}ullstellensatz}, J. Amer. Math.
  Soc. {\bf 1} (1988), no.~4, 963--975.

\bibitem[Lak91]{lak91}
Y.~N. Lakshman, {\em A simple exponential bound on the complexity of computing
  gr{\"o}bner bases of zero-dimensional ideals}, Effective methods in algebraic
  geometry (Castiglioncello, 1990), Progr. Math., vol.~94, Birkhauser Boston,
  1991, pp.~227--234.

\bibitem[Laz87]{laz87}
Robert Lazarsfeld, {\em A sharp {C}astelnuovo bound for smooth surfaces}, Duke
  Math. J. {\bf 55} (1987), 423--429.

\bibitem[{Len}92]{len92}
H.~W. {Lenstra, Jr.}, {\em Algorithms in algebraic number theory}, Bull. Amer.
  Math. Soc. (N.S.) {\bf 26} (1992), no.~2, 211--244.

\bibitem[Lip84]{lip84}
Joseph Lipman, {\em Dualizing sheaves, differentials and residues on algebraic
  varieties}, Ast{\'e}risque, vol. 117, 1984.

\bibitem[LL91]{ll91}
Y.~N. Lakshman and D.~Lazard, {\em On the complexity of zero-dimensional
  algebraic systems}, Effective methods in algebraic geometry (Castiglioncello,
  1990), Progr. Math., vol.~94, Birkhauser Boston, 1991, pp.~217--225.

\bibitem[LLL82]{lll82}
A.~K. Lenstra, H.~W. {Lenstra, Jr.}, and L.~Lov{\'a}sz, {\em Factoring
  polynomials with rational coefficients}, Math. Ann. {\bf 261} (1982),
  515--534.

\bibitem[Mac27]{mac27}
F.~S. Macaulay, {\em Some properties of enumeration in the theory of modular
  systems}, Proc. London Math. Soc. {\bf 26} (1927), 531--555.

\bibitem[MM82]{mm82}
Ernst~W. Mayr and Albert~R. Meyer, {\em The complexity of the word problem for
  commutative semigroups and polynomial ideals}, Adv. in Math. {\bf 46} (1982),
  305--329.

\bibitem[MM83]{mm83}
H.~Michael M{\"o}ller and Ferdinando Mora, {\em Upper and lower bounds for the
  degree of {G}r{\"o}bner bases}, Computer Algebra ({EUROCAL} 83), Lecture
  Notes in Computer Science, vol. 162, Springer-Verlag, 1983, pp.~157--167.

\bibitem[MM86]{mm86}
H.~Michael M{\"o}ller and Ferdinando Mora, {\em New constructive methods in
  classical ideal theory}, J. Algebra {\bf 100} (1986), no.~1, 138--178.

\bibitem[MR88]{mr88}
T.~Mora and L.~Robbiano, {\em The {G}r{\"o}bner fan of an ideal}, J. Symb.
  Comput. {\bf 6} (1988), no.~2--3, 183--208, reprinted in \cite{rob89}.

\bibitem[Mum66]{mum66}
David Mumford, {\em Lectures on curves on an algebraic surface}, Princeton
  University Press, Princeton, New Jersey, 1966.

\bibitem[Mum70a]{mum70a}
David Mumford, {\em Abelian varieties}, Oxford University Press, Oxford, 1970.

\bibitem[Mum70b]{mum70b}
David Mumford, {\em Varieties defined by quadratic equations}, Questions on
  Algebraic Varieties, Centro Internationale Matematica Estivo, Cremonese,
  Rome, 1970, pp.~29--100.

\bibitem[MW83]{mw83}
D.~W. Masser and G.~W{\"u}stholz, {\em Fields of large transcendence degree
  generated by values of elliptic functions}, Invent. Math. {\bf 72} (1983),
  407--464.

\bibitem[Pin86]{pin86}
Henry~C. Pinkham, {\em A {C}astelnuovo bound for smooth surfaces}, Invent.
  Math. {\bf 83} (1986), 491--506.

\bibitem[Ran90]{ran90}
Ziv Ran, {\em Local differential geometry and generic projections of
  threefolds}, J. Differential Geom. {\bf 32} (1990), 131--137.

\bibitem[Rav90]{rav90}
M.~S. Ravi, {\em Regularity of ideals and their radicals}, Manuscripta Math.
  {\bf 68} (1990), 77--87.

\bibitem[Ric74]{ric74}
F.~Richman, {\em Constructive aspects of {N}oetherian rings}, Proc. Amer. Math.
  Soc. {\bf 44} (1974), 436--441.

\bibitem[Rob85]{rob85}
L.~Robbiano, {\em Term orderings on the polynomial ring}, Proceedings of
  EUROCAL '85 (Linz), Lecture Notes in Computer Science, vol. 204,
  Springer-Verlag, 1985, pp.~513--517.

\bibitem[Rob89]{rob89}
Lorenzo Robbiano (ed.), {\em Computational aspects of commutative algebra},
  Academic Press, 1989, {I}SBN 0-12-589590-9.

\bibitem[Sch80]{sch80}
Frank-Olaf Schreyer, {\em {D}ie {B}erechnung von {S}yzygien mit dem
  verallgemeinerten {W}eierstrass'schen {D}ivisionssatz}, Diplomarbeit am
  {F}achbereich {M}athematik der {U}niversit{\"a}t {H}amburg, 1980.

\bibitem[Sch91]{sch91}
Frank-Olaf Schreyer, {\em A standard basis approach to syzygies of canonical
  curves}, J. Reine Angew. Math. {\bf 421} (1991), 83--123.

\bibitem[Ser55]{ser55}
J.-P. Serre, {\em Faisceaux algebrique coherents}, Ann. of Math. (2) {\bf 61}
  (1955), 197--278.

\bibitem[Spe77]{spe77}
D.~Spear, {\em A constructive approach to commutative ring theory}, Proceedings
  of the 1977 {MACSYMA} Users' Conference, NASA CP-2012, 1977, pp.~369--376.

\bibitem[Tei89]{tei89}
Jeremy Teitelbaum, {\em On the computational complexity of the resolution of
  plane curve singularities}, Symbolic and algebraic computation (Rome, 1988),
  Lecture Notes in Computer Science, vol. 358, Springer, 1989, pp.~285--292.

\bibitem[Tei90]{tei90}
Jeremy Teitelbaum, {\em The computational complexity of the resolution of plane
  curve singularities}, Math. Comp. {\bf 54} (1990), no.~190, 797--837.

\bibitem[Tri78]{tri78}
W.~Trinks, {\em {\"U}ber {B}. {B}uchberger's {V}erfahren, {S}ysteme
  algebraischer {G}leichungen zu l{\"o}sen}, J. Number Theory {\bf 10} (1978),
  475--488.

\bibitem[Zac78]{zac78}
G.~Zacharias, Bachelor's thesis, Mass. Inst. of Technology, 1978.

\bibitem[Zar69]{zar69}
Oscar Zariski, {\em An introduction to the theory of algebraic surfaces},
  Lecture Notes in Math., vol.~83, Springer-Verlag, 1969.

\bibitem[ZS76]{zs76}
Oscar Zariski and Pierre Samuel, {\em Commutative algebra, {V}olumes {I},
  {II}}, Graduate texts in mathematics, vol. 28--29, Springer-Verlag,
  1975--1976.

\end{thebibliography}
\newcommand{\noopsort}[1]{} \newcommand{\printfirst}[2]{#1}
  \newcommand{\singleletter}[1]{#1} \newcommand{\switchargs}[2]{#2#1}
\ifx\undefined\bysame
\newcommand{\bysame}{\leavevmode\hbox to3em{\hrulefill}\,}
\fi


\end{document}